\documentclass[12 pt]{amsart}
\title{Derived Equivalences of Twisted Supersingular K3 Surfaces}
\author{Daniel Bragg}
\address{Department of Mathematics, University of California, Berkeley, 
Berkeley, CA 94720}
\email{braggdan@berkeley.edu}
\usepackage{macros}
\usepackage{mathtools}
\begin{document}
\maketitle

\begin{abstract}
 We study the derived categories of twisted supersingular K3 surfaces. We prove a derived crystalline Torelli theorem for twisted supersingular K3 surfaces, characterizing Fourier-Mukai equivalences in terms of the twisted K3 crystals introduced in \cite{BL17}. This is a positive characteristic analog of the Hodge-theoretic derived Torelli theorem of Orlov \cite{MR1465519} and its extension to twisted K3 surfaces by Huybrechts and Stellari \cite{HS04,MR2310257}. We give applications to various questions concerning Fourier-Mukai partners, extending results of C\u{a}ld\u{a}raru \cite{Cal01} and Huybrechts and Stellari \cite{HS04}. We also give an exact formula for the number of twisted Fourier-Mukai partners of a twisted supersingular K3 surface.
\end{abstract}

\tableofcontents


\section{Introduction}

If $X$ is a smooth projective variety, we write $D(X)$ for the bounded derived category of coherent sheaves on $X$. This object is studied as an invariant attached to $X$. The particular case of K3 surfaces over the complex numbers has received extensive attention, and has proven to be both interesting and tractable. A major result is Orlov's Torelli theorem \cite{MR1465519} for the derived category, which characterizes derived equivalences between K3 surfaces over the complex numbers in terms of isomorphisms of certain associated Hodge structures. This result was extended to the case of twisted K3 surfaces by C\u{a}ld\u{a}raru \cite{Cal01} and Huybrechts and Stellari \cite{HS04,MR2310257}. These results allow seemingly difficult questions about derived equivalences to be rephrased in terms of Hodge structures, which are much easier to compute with.

A major obstacle to obtaining similar results for K3 surfaces in positive characteristic is that the classical Torelli theorem no longer applies (this result is a key input to both Orlov's and Huybrechts and Stellari's derived Torelli theorems). A partial replacement was found by Lieblich and Olsson in \cite{LO15}, who introduced a certain filtration on the derived category. They used this to extend results on derived equivalences of K3 surfaces to positive characteristic by lifting to characteristic 0, for instance, the finiteness of the number of Fourier-Mukai partners. 

In this paper, we focus our attention on supersingular K3 surfaces, which are a particularly interesting class of K3 surfaces in positive characteristic. By work of Ogus \cite{Ogus78,Ogus83} and Rudakov and Shafarevich \cite{RTS83}, it is known that supersingular K3 surfaces satisfy a certain Torelli-type result, phrased in terms of crystalline cohomology. The form of this result (reviewed in Section \ref{sec:k3crystals}) is very similar to that of the classical Hodge-theoretic Torelli theorem.


 We use this crystalline Torelli theorem to investigate the derived categories of twisted supersingular K3 surfaces. We are able to obtain similarly strong results as are known for K3 surfaces over the complex numbers, although our methods are necessarily limited to the supersingular case. Our key technical tool is a certain crystalline analog of the twisted Hodge structures of Huybrechts and Stellari \cite{HS04}, introduced in \cite{BL17}. This construction attaches to a twisted K3 surface $(X,\alpha)$ a certain K3 crystal $\tH(X/W,B)$. We show that a derived equivalence of twisted K3 surfaces induces an isomorphism of the corresponding K3 crystals, which preserves certain extra structure. We also show that the converse holds. Our main result is the following crystalline Torelli theorem for the derived category of a twisted supersingular K3 surface.
\begin{Theorem*}
  If $(X,\alpha)$ and $(Y,\beta)$ are twisted supersingular K3 surfaces, then there exists a Fourier-Mukai equivalence $D(X,\alpha)\cong D(Y,\beta)$ if and only if there exists an isomorphism $\tH(X/W,B)\cong\tH(Y/W,B')$ of K3 crystals.
\end{Theorem*}
That is, the twisted K3 crystal exactly captures the relationship of being Fourier-Mukai partners. We prove first a more refined statement (Theorem \ref{thm:Main1}) characterizing which maps of K3 crystals are induced by derived equivalences.

To prove this result, we need to extend various related constructions and results to positive characteristic, such as the existence of certain twisted Chern characters and moduli spaces of twisted sheaves. Much of this has been already taken care of in \cite{BL17}. We review the needed definitions and results below. In two appendices, we show that our twisted Chern characters take integral values, and that derived equivalences of certain twisted K3 surfaces in positive characteristic are orientation preserving, extending results of Huybrechts, Macr\`\i, and Stellari \cite{MR2553878} and Reinecke \cite{2017arXiv171100846R}. We also include a comparison of our twisted Chern characters with those of Huybrechts and Stellari \cite{HS04}.

Our main theorem allows the translation of questions about derived equivalences into questions about isomorphisms of supersingular K3 crystals. By results of Ogus (reviewed in Section \ref{sec:characteristicsubspaces}), supersingular K3 crystals are essentially objects of semi-linear algebra, and hence are quite easy to compute with. Using these techniques, we apply our main theorem in Sections \ref{sec:characteristicsubspaces}, \ref{sec:applications}, and \ref{sec:counting} to prove various facts about derived equivalences of twisted supersingular K3 surfaces. For instance, we recover the result of Lieblich and Olsson \cite{LO15} that (non-twisted) supersingular K3 surfaces are derived equivalent if and only if they are isomorphic (see Theorem \ref{thm:Untwisted1}). We also prove that if $(X,\alpha)$ is a twisted supersingular K3 surface then $D(X,\alpha)\cong D(X,\alpha^k)$ for any integer $k$ coprime to the order of $\alpha$ (Theorem \ref{thm:caldararusthing}). Finally, in sections \ref{sec:orthogonal} and \ref{sec:counting} we use characteristic subspaces to count twisted Fourier-Mukai partners. The final result is Theorem \ref{Count2}, which gives the following formula for the number of twisted partners of a twisted supersingular K3 surface $(X,\alpha)$ in terms of two discrete invariants $\sigma_0,m$ attached to $(X,\alpha)$.
  \[
    \#\FM(X,\alpha)=
      \begin{dcases}
        \frac{p^{\sigma_0}+1}{p^m+1}(p^{\sigma_0-1}-1)+1 &\mbox{ if }1\leq\sigma_0\leq 10\\
        \frac{p^{\sigma_0}+1}{p^m+1}(p^{\sigma_0-1}-1)   &\mbox{ if }\sigma_0=11\\
      \end{dcases}
  \]
We remark that Lieblich and Olsson \cite{LO15.2} have also used the crystalline Torelli thoerem to study the derived categories of K3 surfaces in positive characteristic. They obtain results for general K3 surfaces in positive characteristic by deformation to the supersingular locus, replacing the lifting arguments of \cite{LO15}.

\subsection{Notation}
We fix throughout an algebraically closed field $k$ of characteristic $p>0$. Unless otherwise stated, we will assume that $p\neq 2$. If $X$ is a variety over $k$, the notation $\H^i(X,\m_{p^n})$ will always be used to denote the cohomology of $\m_{p^n}$ with respect to the flat (fppf) topology. A \textit{K3 surface} is a smooth projective surface $X$ with trivial canonical bundle and $\H^1(X,\mathscr{O}_X)=0$.

\subsection{Acknowledgements}
The results of this paper were obtained during the course of the author's PhD work at the University of Washington. The author thanks Nicolas Addington, Ben Antieau, Max Lieblich, Martin Olsson, and Emanuel Reinecke for useful conversations and correspondences.

\section{Supersingular K3 surfaces}

Let $X$ be a K3 surface over $k$. We say that $X$ is \textit{supersingular} if $X$ has Picard number 22.\footnote{This is the strongest possible notion of supersingularity, and is sometimes known as \textit{Shioda supersingular}. This condition is equivalent to all of the various weaker notions appearing in the literature (see eg. Chapter 18 of \cite{Huy06}) by work of Charles, Kim, Madapusi Pera, and Maulik \cite{Charles13,Charles14,KMP16, MP15,Mau14}.} There is a beautiful theory of such surfaces, with many special results and techniques. Their systematic study was initiated by Artin \cite{Artin74}. It is known that the discriminant of the N\'{e}ron-Severi group of a supersingular K3 surface is of the form $-p^{2\sigma_0}$ for some integer $1\leq\sigma_0\leq 10$, known as the \textit{Artin invariant} of $X$. The moduli space of supersingular K3 surfaces is $9$ dimensional, and the locus of surfaces $X$ such that $\sigma_0(X)\leq\sigma_0$ is closed of dimension $\sigma_0-1$. In particular, the locus of supersingular K3 surfaces with $\sigma_0=1$ is discrete, and in fact there is up to isomorphism a unique such surface.

\subsection{The Brauer group}
We briefly describe the Brauer group of a supersingular K3 surface. If $X$ is a K3 surface over the complex numbers, or more generally over an algebraically closed field of characteristic 0, then there is an abstract isomorphism of groups
\[
  \Br(X)\cong\left(\bQ/\bZ\right)^{\oplus 22-\rho(X)}
\]
(see for instance Example 18.1.7ff of \cite{Huy16}). In particular, the Brauer group is divisible and countable. Suppose that $X$ is a K3 surface over an algebraically closed field $k$ of positive characteristic. If $X$ has finite height, then $\Br(X)$ is again divisible and countable, and has a similar form. If $X$ is supersingular, then the situation is completely different. Artin showed in \cite{Artin74} that there is an abstract isomorphism of groups
\[
  \Br(X)\cong k
\]
between the Brauer group and the underlying additive group of the ground field $k$. In particular, $\Br(X)$ is not divisible, is $p$-torsion, and has the same cardinality as $k$. The explanation for this distinctive behavior comes from flat cohomology. The flat Kummer sequence gives a short exact sequence
\[
  0\to\Pic(X)\otimes\mathbb{Z}/p\mathbb{Z}\to\H^2(X,\m_p)\to\Br(X)[p]\to 0
\]
Artin showed that $\H^2(X,\m_p)$ is the $k$-points of a certain naturally associated algebraic group, and defined the subgroup $\U^2(X,\m_p)\subset\H^2(X,\m_p)$ to be the $k$-points of the connected component of the identity (we record an alternative characterization of the subgroup $\U^2(X,\m_p)$ in Definition \ref{def:alternatedef}). The completion of this algebraic group at the identity is identified with the formal Brauer group of $X$, which is isomorphic to $\widehat{\bG}_a$. This implies that we have an isomorphism $\U^2(X,\m_p)\cong\bG_a(k)$. Let $N(X)=\Pic(X)$ be the N\'{e}ron-Severi group of $X$. Artin showed that the Kummer sequence induces a diagram
\begin{equation}\label{eq:artinssubgroup}
  \begin{tikzcd}
    0\arrow{r}&\dfrac{pN(X)^{\vee}}{pN(X)}\arrow{r}\arrow[hook]{d}&\U^2(X,\m_p)\arrow{r}\arrow[hook]{d}&\Br(X)\isor{d}{}\arrow{r}&0\\
    0\arrow{r}&\dfrac{N(X)}{pN(X)}\arrow{r}&\H^2(X,\m_p)\arrow{r}&\Br(X)\arrow{r}&0
  \end{tikzcd}
\end{equation}
In particular, the map $\U^2(X,\m_p)\to\Br(X)$ is surjective and has finite kernel, explaining the profusion of $p$-torsion Brauer classes in the supersingular case.

Despite this very non-classical behavior of the Brauer group, we will see that many structures associated to twisted K3 surfaces over the complex numbers have analogs for supersingular K3 surfaces.

\subsection{K3 crystals and Ogus's crystalline Torelli theorem}\label{sec:k3crystals}

Let $W$ denote the ring of Witt vectors of $k$, let $K=W[p^{-1}]$ be its field of fractions, and let $\sigma\colon W\to W$ be its Frobenius. If $X$ is any K3 surface over $k$ (not necessarily supersingular), then the second crystalline cohomology group $\H^2(X/W)$ of $X$ is a free $W$-module of rank 22. It is equipped with a perfect pairing
\begin{equation}\label{eq:cup product}
  \H^2(X/W)\otimes_W\H^2(X/W)\to W
\end{equation}
and an endomorphism $\Phi\colon \H^2(X/W)\to \H^2(X/W)$ which is $\sigma$-linear,\footnote{That is, $\Phi(\lambda h)=\sigma(\lambda)\Phi(h)$ for all $h\in\H^2(X/W)$ and $\lambda\in W$.} and satisfies a certain compatibility with the pairing. This object may be viewed as a positive characteristic analog of the singular cohomology group $\H^2(X,\bZ)$ of a K3 surface over the complex numbers, equipped with its pairing and Hodge structure. Following Ogus (see \cite[Definition 1]{Ogus83}), we abstract the above structure in the following definition.
\begin{Definition}\label{def:K3crystal}
A \textit{K3 crystal of rank $n$} is a $W$-module $H$ of rank $n$ equipped with a $\sigma$-linear endomorphism $\Phi\colon H\to H$ and a symmetric bilinear form $H\otimes_WH\to W$ such that
  \begin{enumerate}
      \item $p^2H\subset\Phi(H)$,
      \item $\Phi\otimes k$ has rank 1,
      \item the pairing $\langle\_,\_\rangle$ is perfect,
      \item $\langle\Phi(x),\Phi(y)\rangle=p^2\sigma\langle x,y\rangle$ for all $x,y\in H$, and
      \item the crystalline discriminant of $(H,\Phi)$ is $-1$.
  \end{enumerate}
 We will use the notation $\varphi=p^{-1}\Phi$.
\end{Definition}

The second crystalline cohomology group $\H^2(X/W)$ of a K3 surface equipped with the cup product pairing and the endomorphism $\Phi$ satisfies the above properties, and is the basic example of a K3 crystal (of rank 22) (see for instance \cite{Ogus83}).
\begin{Definition}
  The \textit{Tate module} of a K3 crystal $H$ is
  \[
    T_H=\left\{h\in H|\Phi(h)=ph\right\}
  \]
  This has a natural structure of $\bZ_p$-module, and is equipped with the restriction of the bilinear form on $H$.
\end{Definition}
We recall some abstract structural results. An $F$-\textit{crystal} is a pair $(H,\Phi)$ where $H$ is a free $W$-module of finite rank and $\Phi\colon H\to H$ is a $\sigma$-linear endomorphism. A map $(H,\Phi)\to (H',\Phi')$ of $F$-crystals is a map $H\to H'$ of $W$-modules that commutes with the respective semilinear operators. Such a map is said to be an \textit{isogeny} if it induces an isomorphism $H\otimes K\to H'\otimes K$. There is a well known classification of $F$-crystals up to isogeny in terms of slopes, due to Dieudonn\'{e}-Manin. Note that any K3 crystal gives rise in particular to an $F$-crystal by forgetting the pairing. The isogeny type of this crystal is constrained by the extra structure imposed above, and gives rise to an invariant of the K3 crystal. We will mostly be interested in the following case in this paper.
\begin{Definition}
  A K3 crystal $H$ is \textit{supersingular} if it has pure slope 1.
\end{Definition}
Equivalently, a K3 crystal $H$ is supersingular exactly when the inclusion $T_H\otimes W\to H$ is an isogeny (that is, the map $T_H\otimes K\to H\otimes K$ is an isomorphism). If $X$ is a K3 surface, then the K3 crystal $\H^2(X/W)$ is supersingular if and only if $X$ is supersingular.

If $H$ is supersingular, then by Proposition 3.13 of \cite{Ogus78} the discriminant of the pairing on $T_H$ is equal to $-p^{2\sigma_0}$ for some integer $\sigma_0\geq 1$.
\begin{Definition}\label{def:artininvariant1}
  The \textit{Artin invariant} $\sigma_0(H)$ of a supersingular K3 crystal $H$ is the integer $\sigma_0$.
\end{Definition}
It is immediate that if $X$ is a supersingular K3 surface, then the Artin invariant of $X$ is equal to the Artin invariant of $\H^2(X/W)$. Suppose that $X$ and $Y$ are supersingular K3 surfaces. An obvious necessary condition for a map $\H^2(X/W)\iso\H^2(Y/W)$ of K3 crystals to be induced by a morphism is that it maps an ample class to an ample class. It is easy to see that this need not always be the case. In fact, it need not even be the case that the Picard group of $X$ is mapped to the Picard group of $Y$ (unlike in the complex case, where this is automatic). We will therefore need to track the inclusion of the Picard group as a piece of extra data.

\begin{Definition}\label{def:marking}
  A \textit{marked K3 crystal} is a pair $(H,N)$ where $H$ is a K3 crystal and $N$ is a lattice equipped with an isometric inclusion $N\to H$ that induces an isomorphism $N\otimes\bZ_p\iso T_H$.
  
  An isomorphism of marked K3 crystals $(H,N)\to(H',N')$ consists of an isomorphism $H\to H'$ of K3 crystals that restricts to an isomorphism $N\to N'$.
\end{Definition}
The first Chern class map $\pic(X)\to\H^2(X/W)$ gives a marking of $\H^2(X/W)$. The following is Ogus's crystalline Torelli theorem.\footnote{Ogus proved this result under the assumption that $p\geq 5$. It has been extended to $p\geq 3$ in \cite{BL17}}
\begin{Theorem}\label{thm:Ogus1}\cite[Theorem II]{Ogus83}
  Let $X$ and $Y$ be supersingular K3 surfaces. If 
  \[
    f\colon (\H^2(X/W),\pic(X))\iso (\H^2(Y/W),\pic(Y))
  \]
  is an isomorphism of marked supersingular K3 crystals such that the map $\pic(X)\to\pic(Y)$ sends an ample class to an ample class, then $f$ is induced by a unique isomorphism $X\iso Y$.
\end{Theorem}

Somewhat surprisingly, it turns out that if we are only interested in the existence of an isomorphism, we may neglect the marking entirely.

\begin{Theorem}\label{thm:Ogus2}\cite[Theorem I]{Ogus83}
  If $X$ and $Y$ are supersingular K3 surfaces, then $X$ is isomorphic to $Y$ if and only if there exists an isomorphism $\H^2(X/W)\iso\H^2(Y/W)$ of K3 crystals.
\end{Theorem}

\section{Twisted K3 crystals and derived equivalences}

In the following sections we present the constructions and definitions needed to state our main theorem. 
\subsection{Crystalline B-fields}

We review the construction of crystalline B-fields defined in Section 3.4 of \cite{BL17}. Let $X$ be a K3 surface over $k$. We begin by relating the flat cohomology of $\m_{p}$ to the de Rham cohomology of $X$. Consider the Kummer sequence
\[
  1\to\m_{p}\to\bG_m\xrightarrow{x\mapsto x^{p}}\bG_m\to 1
\]
which is exact in the flat topology. Let $\varepsilon\colon X_{\fl}\to X_{\etale}$ be the natural map from the big fppf site of $X$ to the small \'{e}tale site of $X$. By a theorem of Grothendieck, the cohomology of the complex $\R\varepsilon_*\bG_m$ vanishes in all positive degrees. Because $X$ is smooth over a perfect scheme, the restriction of $\m_{p}$ to the small \'{e}tale site of $X$ is trivial. Applying $\varepsilon_*$ to the Kummer sequence, we obtain an exact sequence
\[
  1\to\bG_m\xrightarrow{x\mapsto x^{p}}\bG_m\to\R^1\varepsilon_*\m_{p}\to 1
\]
of sheaves on the small \'{e}tale site of $X$. It follows that
\[
  \R^1\varepsilon_*\m_{p}=\bG_m/\bG_m^{p}
\]
where the quotient is taken in the \'{e}tale topology. We therefore obtain isomorphisms
\begin{equation}\label{eq:flattoetale}
  \H^i(X_{\fl},\m_{p})\iso\H^{i-1}(X_{\etale},\bG_m/\bG_m^{p})
\end{equation}
We consider the map of \'{e}tale sheaves
\begin{equation}\label{eq:dlog}
  d\log\colon\bG_m\to\Omega^1_X
\end{equation}
defined by $f\mapsto df/f$. This map sends any $p$-th power to 0, and hence descends to the quotient. Combined with the isomorphism (\ref{eq:flattoetale}), we find a map
\begin{equation}\label{eq:dlog2}
  d\log\colon\H^2(X,\m_p)\hto\H^2_{dR}(X/k)
\end{equation}
(the injectivity follows for instance from Corollaire 2.1.18 of \cite{Ill79}). The essential property of this map that we will need is that it is compatible with the de Rham Chern character, in the sense that the diagram
\[
  \begin{tikzcd}
    \H^1(X,\G_m)\arrow{dr}{c_1^{dR}}\arrow{d}&\\
    \H^2(X,\m_p)\arrow[hook]{r}{d\log}&\H^2_{dR}(X/k)
  \end{tikzcd}
\]
commutes, where the vertical arrow is the boundary map induced by the Kummer sequence. As the third crystalline cohomology group of $X$ is torsion free, the canonical map
\[
  \H^2(X/W)\otimes k\to\H^2_{dR}(X/k)
\]
is an isomorphism. Combined with (\ref{eq:dlog2}), we arrive at a diagram
\begin{equation}\label{eq:crystallinediagram}
  \begin{tikzcd}[row sep=small]
    \H^2(X/W)\arrow[two heads]{r}{\pi}&\H^2_{dR}(X/k)&\\
    &\H^2(X,\m_{p})\arrow[phantom]{u}{\rotatebox{90}{$\subset$}}\arrow[two heads]{r}&\Br(X)[p]
  \end{tikzcd}
\end{equation}
Write $\H^2(X/K)=\H^2(X/W)\otimes_WK$.
\begin{Definition}
  Let $\alpha\in\Br(X)[p]$ be a $p$-torsion Brauer class. A \textit{B-field lift} of $\alpha$ is an element $B=\frac{a}{p}\in\H^2(X/K)$, where $a$ is an element of $\H^2(X/W)$ whose image in $\H^2_{dR}(X/k)$ is contained in the image of $d\log$, and maps to $\alpha$ under the lower horizontal map of (\ref{eq:crystallinediagram}).
\end{Definition}
We note that, unlike in the Hodge-theoretic setting, it is no longer true that every class in the rational crystalline cohomology group is a B-field lift of some Brauer class. We give in \cite{BL17} the following characterization of B-fields.
\begin{Lemma}\cite[Lemma 3.4.11]{BL17}\label{lem:Bfields!}
  A class $B\in p^{-1}\H^2(X/W)$ is a B-field lift of a class in $\Br(X)[p]$ if and only if
  \[
    B-\varphi(B)\in\H^2(X/W)+\varphi(\H^2(X/W))
  \]
  (recall that $\varphi=p^{-1}\Phi$)
\end{Lemma}
\begin{Remark}
  It is possible to make a similar definition of a B-field lift associated to a class in $\Br(X)[p^n]$ by replacing de Rham cohomology with de Rham-Witt cohomology. One can show that a similar characterization to that of Lemma \ref{lem:Bfields!} holds. As the Brauer group of a supersingular K3 surface is $p$-torsion, the above definition will suffice for our applications in this paper.
\end{Remark}

\subsection{The twisted Mukai crystal}

In this section we will show how to associate to a pair $(X,\alpha)$, where $X$ is a K3 surface and $\alpha\in\Br(X)[p]$, a certain K3 crystal of rank 24. The material in this section is essentially contained in Section 3.4 of \cite{BL17}, although with slightly different notation. We begin by defining the Mukai crystal associated to a K3 surface $X$, as originally introduced in \cite{LO15}.

\begin{Definition}\label{def:MukaiCrystal}
The \textit{Mukai crystal} of $X$ is the $W$-module
\[
  \tH(X/W)=\H^0(X/W)(\text{-}1)\oplus \H^2(X/W)\oplus \H^4(X/W)(1)
\]
equipped with the twisted Frobenius $\tPhi:\tH(X/W)\to\tH(X/W)$. We define the \textit{Mukai pairing} on $\tH(X/W)$ by
\[
  \langle(a,b,c),(a',b',c')\rangle=-ac'+b.b'-a'c
\]
where $b.b'$ is the cup product pairing~\eqref{eq:cup product}.
\end{Definition}
Both $\H^0(X/W)$ and $\H^4(X/W)$ are canonically isomorphic (as $W$-modules) to $W$. Under these identifications, the twisted Frobenius action is given by
\[
  \tPhi(a,b,c)=(p\sigma(a),\Phi(b),p\sigma(c))
\]
It follows immediately from the definitions that $\tH(X/W)$ is a K3 crystal of rank 24.

We next explain how to twist the crystal $\tH(X/W)$ by a B-field. Given any element $B\in\H^2(X/K)$ (for instance, a B-field), cupping with the class $e^B=(1,B,B^2/2)$ defines an isometry
\[
   e^B\colon \tH(X/K)\to\tH(X/K)
\]
where $\tH(X/K)=\tH(X/W)\otimes_WK$. Explicitly, this map is given by
\begin{equation}\label{eq:cupping with B}
  e^B(a,b,c)=\left(a,b+aB,c+b.B+a\frac{B^2}{2}\right)
\end{equation}
\begin{Definition}\label{def:twistedK3crystal}
Let $\alpha$ be a class in $\Br(X)[p]$ and let $B$ be a B-field lift of $\alpha$. The \textit{twisted Mukai crystal} associated to $(X,\alpha)$ and $B$ is
\[
  \tH(X/W,B)=e^B\tH(X/W)\subset \tH(X/K)
\]
\end{Definition}

See Remark \ref{rem:conventions and HS} for the relation between this definition and the twisted Mukai lattice of Huybrechts and Stellari. It is shown in \cite{BL17} that $\tH(X/W,B)$ has a natural structure of K3 crystal.
\begin{Theorem}\cite[Proposition 3.4.15]{BL17}\label{thm:itsacrystal}
  The submodule $\tH(X/W,B)\subset\tH(X/K)$ is preserved by the Frobenius endomorphism $\tPhi$. The $W$-module $\tH(X/W,B)$ equipped with the endomorphism $\tPhi$ and the restriction of the Mukai pairing on $\tH(X/K)$ is a K3 crystal of rank 24, and is supersingular if and only if $X$ is supersingular.
\end{Theorem}
The most difficult condition to verify is that $\tH(X/W,B)$ is preserved by $\tPhi$, which follows from the characterization of B-fields given in Lemma \ref{lem:Bfields!}. 

Although the embedding $\tH(X/W,B)\subset\tH(X/K)$ does depend on the choice of B-field, we will show that the isomorphism class of the crystal $\tH(X/W,B)$ is independent of the choice of B-field.
\begin{Lemma}\label{lem:independent}
  If $B$ and $B'$ are two different B-field lifts of $\alpha$, then there exists an isomorphism of K3 crystals $\tH(X/W,B)\cong\tH(X/W,B')$.
\end{Lemma}
\begin{proof}
While cupping with $e^{B'-B}$ defines an isometry between $\tH(X/W,B)$ and $\tH(X/W,B')$, this map does not in general commute with the Frobenius operators, and so does not define an isomorphism of K3 crystals. Indeed, a direct computation yields the useful identity
\begin{equation}\label{eq:identity1}
  \tPhi(e^A(a,b,c))=e^{\varphi(A)}\tPhi (a,b,c)
\end{equation}
valid for any $A\in\H^2(X/K)$ and $(a,b,c)\in\tH(X/W)$ (recall that $\varphi=p^{-1}\Phi$). Instead, consider the diagram
\[
  \begin{tikzcd}[row sep=small]
    &&\H^2_{dR}(X/k)&&\\
    0\arrow{r}&\Pic(X)\otimes\bZ/p\bZ\arrow{ur}{c_1^{dR}}\arrow{r}&\H^2(X,\m_p)\arrow[phantom]{u}{\rotatebox{90}{$\subset$}}\arrow{r}&\Br(X)[p]\arrow{r}&0
  \end{tikzcd}
\]
where the lower row is exact. The image of $pB'-pB$ in $\H^2_{dR}(X/k)$ maps to 0 in $\Br(X)[p]$, and hence is equal to $c_1^{dR}(\sL)$ for some line bundle $\sL$. Thus, writing $t=c_1^{\cry}(\sL)$, we have
\[
  B'-B=\frac{t}{p}+h
\]
for some $h\in\H^2(X/W)$. Note that if $h\in \H^2(X/W)$ then $e^h\in \H^*(X/W)$. Now, $\varphi(t/p)=t/p$, so by (\ref{eq:identity1}) cupping with $e^{t/p}$ defines an isomorphism
\[
  e^{t/p}\colon \tH(X/W,B)\iso e^{B'-h}\tH(X/W)=\tH(X/W,B')
\]
of K3 crystals.
\end{proof}
\begin{Remark}
  As a submodule of $\tH(X/K)$, the twisted crystal $\tH(X/W,B)$ depends only on the image of $B$ in $\H^2(X,\m_p)$. In other words, $\tH(X/W,B)$ may be viewed as a structure canonically attached to a class in $\H^2(X,\m_p)$. This is the viewpoint taken in \cite{BL17}.
\end{Remark}
We can use the crystal $\tH(X/W,B)$ to extend the Artin invariant to twisted supersingular K3 surfaces.
\begin{Definition}
  The \textit{Artin invariant} of a twisted supersingular K3 surface $(X,\alpha)$ is the Artin invariant of the supersingular K3 crystal $\tH(X/W,B)$ (in the sense of Definition \ref{def:artininvariant1}).
\end{Definition}
In Section 3.4 of \cite{BL17}, we make the following computation.
\begin{Proposition}\cite[Corollary 3.4.23]{BL17}
  If $(X,\alpha)$ is a twisted supersingular K3 surface, then
  \[
    \sigma_0(X,\alpha)=\begin{cases}
    \sigma_0(X)+1&\mbox{ if }\alpha\neq 0\\
    \sigma_0(X)  &\mbox{ if }\alpha=0
    \end{cases}
  \]
\end{Proposition}
In particular, the Artin invariant of a twisted supersingular K3 surface is an integer $1\leq\sigma_0\leq 11$.

\subsection{Twisted Chern characters and the twisted N\'{e}ron-Severi group}

Let $X$ be a K3 surface over an algebraically closed field $k$ of arbitrary characteristic. Let $N(X)$ be the N\'{e}ron-Severi group of $X$ (recall that if $X$ is a K3 surface then $N(X)=\Pic(X)$).
\begin{Definition}
  We define the \textit{extended N\'{e}ron-Severi lattice} of $X$ by
  \[
    \tN(X)=\langle (1,0,0)\rangle\oplus N(X)\oplus\langle (0,0,1)\rangle
  \]
  which we equip with the Mukai pairing
  \[
    \tN(X)\otimes\tN(X)\to\bZ
  \]
  given by the formula $\langle (a,b,c),(a',b',c')\rangle=-a.c'+b.b'-a'.c$.
\end{Definition}
We define the \textit{Mukai vector} of a coherent sheaf $\sE$ on $X$ by
\[
  v(\sE)=\ch(\sE).\sqrt{\Td(X)}
\]
Recall that we have $\sqrt{\Td(X)}=(1,0,1)$, and that the Chern characters of a coherent sheaf on a K3 surface are integral. The Mukai vector thus gives rise to a map
\[
  v(\_):K(X)\to\tN(X)
\]
which is an isomorphism of groups. If we equip $K(X)$ with the pairing induced by
\begin{equation}\label{eq:pairing}
  \langle [\sE],[\sF]\rangle=-\sum_i(-1)^i\dim_k\Ext^i(\sE,\sF)
\end{equation}
it is even an isometry.

We wish to replicate these structures in the twisted setting. Let $X$ be a smooth projective variety and $\alpha\in\Br(X)$ a Brauer class (we have in mind the case where $X$ is a K3 surface, or the product of two K3 surfaces). Before even considering the category of $\alpha$-twisted sheaves, we must make a non-canonical choice, either of a cocycle or of a $\mathbb{G}_m$-gerbe representing the cohomology class $\alpha$. We will opt for the latter, although either would be sufficient for our purposes. A good reference for the gerby perspective is \cite[12.3]{MR3495343}. For the readers convenience, we briefly discuss the cocycle perspective in Appendix A. 

There are several different definitions of Chern characters for twisted sheaves appearing in the literature. These are all essentially equivalent, but differ by various shifts and signs which can be confusing to compare. The particular choice which we will employ in the present document is essentially that used in \cite{BL17,LMS11}. We recall the definition below. In Appendix A, we include a comparison with the twisted Chern characters used by Huybrechts and Stellari \cite{HS04}.

 Let $\pi\colon\sX\to X$ be a $\bG_m$-gerbe whose associated cohomology class is $\alpha$. If $n\alpha=0$, then there exists an $n$-fold twisted invertible sheaf, say $\sL$, on $\sX$. Choosing such a sheaf allows us to compare $n$-fold twisted sheaves on $\sX$ and sheaves on $X$.
\begin{Definition}\label{def:twistedcherncharacter}
  Let $\sE$ be a locally free twisted sheaf of positive rank on $\sX$. The \textit{twisted Chern character} of $\sE$ (with respect to $\sL$) is
  \[
    \ch^{\sL}(\sE)=\sqrt[\leftroot{-2}\uproot{2}n]{\ch(\pi_*(\sE^{\otimes n}\otimes\sL^{\vee}))}\in A^*(X)\otimes\bQ
  \]
    where by convention we choose the $n$-th root so that $\rk(\ch^{\sL}(\sE))=\rk(\sE)$. We define the twisted Mukai vector of $\sE$ by
  \[
    v^{\sL}(\sE)=\ch^{\sL}(\sE).\sqrt{\Td(X)}
  \]
\end{Definition}
Any twisted sheaf on $\sX$ admits a resolution by locally free twisted sheaves of positive rank, so the above definitions of the twisted Chern character and twisted Mukai vector extend by additivity to all of $K(X,\alpha)$.

The choice of an $n$-fold twisted invertible sheaf $\sL$ on $\sX$ is essentially the same as a choice of preimage of $\alpha$ under the map
\begin{equation}\label{eq:map1}
  \H^2(X,\m_n)\twoheadrightarrow\Br(X)[n],
\end{equation}
Indeed, giving an $n$-fold twisted invertible sheaf $\sL$ on $\sX$ is the same as giving an isomorphism $\sX^{\wedge n}\iso B\bG_m$ of $\bG_m$-gerbes, which is the same as giving a $\m_n$-gerbe $\sX'$ together with an identification of its associated $\bG_m$-gerbe with $\sX$. The cohomology class of such a gerbe in $\H^2(X,\m_n)$ then maps to $\alpha$ under (\ref{eq:map1}). Conversely, given a preimage $\alpha'$ of $\alpha$, choosing a representing $\m_n$-gerbe and an isomorphism between its associated $\bG_m$-gerbe and $\sX$ gives rise to an invertible $n$-fold twisted sheaf $\sL$ on $\sX$. Moreover, the resulting twisted Chern character $\ch^{\sL}$ depends only on $\alpha'$.
\begin{Definition}\label{def:twisted Chern character part 2}
  If $\alpha'$ is a preimage of $\alpha$ under (\ref{eq:map1}), we write $\ch^{\alpha'}$ for the twisted Chern character determined by $\alpha'$, and set $v^{\alpha'}=\ch^{\alpha'}.\sqrt{\Td(X)}$.
\end{Definition}

Let us suppose now that the characteristic of $k$ is $p$ and that $n=p$. Let $B$ be a crystalline B-field lift of $\alpha$. By definition, the reduction modulo $p$ of $pB$ is an element $\alpha'\in\H^2(X,\m_p)$ which maps to $\alpha$, and we write $\ch^{B}=\ch^{\alpha'}$ (resp. $v^B=v^{\alpha'}$) for the resulting twisted Chern character (resp. twisted Mukai vector). Similarly, if $X\times Y$ is a product of two K3 surfaces, $\alpha\in\Br(X)$ and $\beta\in\Br(Y)$ are Brauer classes or order $p$ with B-field lifts $B,B'$, we define $\ch^{-B\boxplus B'}$ and $v^{-B\boxplus B'}$.

We will identify a natural subgroup of $A^*(X)\otimes\bQ$ which will be a recipient for $\ch^B$ and $v^B$, and so provide a replacement for the twisted N\'{e}ron-Severi group $\tN(X)$.
\begin{Definition}
  Let $X$ be a K3 surface, $\alpha\in\Br(X)[p]$ a $p$-torsion Brauer class, and $B$ a B-field lift of $\alpha$. The \textit{twisted N\'{e}ron-Severi lattice} associated to $(X,\alpha,B)$ is
  \[
    \tN(X,B)=(\tN(X)\otimes\bZ[p^{-1}])\cap \exp(B)\tH(X/W)
  \]
  We equip $\tN(X,B)$ with the restriction of the Mukai pairing on $\tH(X/W,B)$.
\end{Definition}

We mention two results concerning this object.
\begin{Proposition}\label{prop:integrality1}
  If $\sE$ is an $\alpha$-twisted sheaf, then $\ch^B(\sE)$ and $v^B(\sE)$ are contained in $\tN(X,B)$. The resulting map
  \[
    v^B(\_)\colon K(X,\alpha)\to\tN(X,B)
  \]
  is an isomorphism of groups. If we equip the left hand side with the pairing (\ref{eq:pairing}), it is an isometry.
\end{Proposition}
We defer the proof to Appendix \ref{app:twistedChernCharacters}. See also Proposition 4.1.9 of \cite{BL17} for a different proof.

\begin{Remark}\label{rem:conventions and HS}
    The object analogous to $\tN(X,B)$ in the theory of Huybrechts and Stellari is the \textit{twisted Picard group}, defined on page 12 of \cite{HS04}. With our definition of $\tN(X)$, this is given by
    \[
        \Pic(X,B)=(\exp(B)\tN(X)\otimes\mathbb{Q})\cap\tH(X,\mathbb{Z})
    \]
   Due to differing conventions, this definition differs from ours in two separate ways. First, there is an overall shift by a factor of $\exp (-B)$. This is because we have chosen to define our lattice $\tH(X/W,B)$ as a sublattice of the rational cohomology $\tH(X/K)$ (in general different from the untwisted Mukai lattice $\tH(X/W)\subset\tH(X/K)$), whereas Huybrechts and Stellari define $\tH(X,B,\mathbb{Z})$ to be equal to the untwisted Mukai lattice $\tH(X,\mathbb{Z})$, but with Hodge structure twisted by $B$. Second, even after shifting by $\exp(-B)$, our definitions differ by replacing $B$ with $-B$. This difference is forced by our use of different conventions for the normalization of Chern characters, as described in Appendix \ref{ssec:comparison}.
\end{Remark}

Let $\tT(X,B)$ denote the Tate module of $\tH(X/W,B)$.
\begin{Proposition}\cite[Proposition 4.1.10]{BL17}
  The inclusion $\tN(X,B)\hto\tH(X/W,B)$ induces an isomorphism
  \[
    \tN(X,B)\otimes\bZ_p\iso\tT(X,B)
  \]
  of $\mathbb{Z}_p$-modules which is an isometry with respect to the natural pairings.
\end{Proposition}
Thus, the inclusion $\tN(X,B)\hto\tH(X/W,B)$ gives a marking of $\tH(X/W,B)$ in the sense of Definition \ref{def:marking}.

\subsection{Action on cohomology}

Let $(X,\alpha)$ and $(Y,\beta)$ be twisted supersingular K3 surfaces, where $\alpha$ has order $n$ and $\beta$ has order $m$. Consider a complex $P^\bullet\in D(X\times Y,-\alpha\boxtimes\beta)$ of twisted sheaves such that the induced Fourier-Mukai functor
\[
  \Phi_{P^\bullet}\colon D(X,\alpha)\to D(Y,\beta)
\]
is an isomorphism. To study such objects, we will consider their action on the twisted Grothendieck groups and the twisted Mukai crystal. Choose B-field lifts $B$ of $\alpha$ and $B'$ of $\beta$. Using the standard formalism, we obtain maps
\begin{align*}
    \Phi^{K}_{\left[P^\bullet\right]}(\_)&\colon K(X,\alpha)\to K(Y,\beta)\\
    \Phi^{\cry}_{v^{-B\boxplus B'}(P^\bullet)}(\_)&\colon \tH(X/K)\to\tH(Y/K)
\end{align*}
These maps are compatible, in the sense that the diagram
\[
    \begin{tikzcd}[column sep=large]
      K(X,\alpha)\arrow{r}{\Phi^K_{[P^\bullet]}}\arrow{d}[swap]{v^B}&K(Y,\beta)\arrow{d}{v^{B'}}\\
      \tH(X/K)\arrow{r}{\Phi^{\cry}_{v^{-B\boxplus B'}(P^{\bullet})}}&\tH(Y/K)
    \end{tikzcd}
\]
commutes (This is shown in \cite[Lemma 4.1.6]{BL17}. See also \cite[Corollary 5.29]{Huy06} for a proof in the untwisted case). Furthermore, the map $\Phi^{\cry}_{v^{-B\boxplus B'}(P^{\bullet})}$ is an isometry with respect to the Mukai pairing and, and commutes with the respective Frobenius operators. To see the former, note that as~\eqref{eq:cupping with B} is an isometry, we may assume that the twists are trivial. In this case, the result follows exactly as in the classical case over the complex numbers (see eg. \cite[Corollary 10.7, Proposition 5.44]{Huy06}). To see the commutativity with the Frobenius operators, we observe that the class $v^{-B\boxplus B'}(P^{\bullet})$ is contained in the subspace $\oplus_i\H^{2i}(X\times Y/K)^{\Phi=p^i}\subset\H^*(X\times Y/K)$. The claim then follows upon expanding in terms of the Kunneth decomposition (this reasoning is identical to that of \cite[Proposition 5.39]{Huy06}, which is the analogous claim over the complex numbers). As an immediate corollary, we obtain the following.
\begin{Proposition}\label{prop:FMpartner1}
  Let $(X,\alpha)$ be a twisted supersingular K3 surface. If $Y$ is a smooth projective variety and there exists a derived equivalence $(X,\alpha)\cong (Y,\beta)$ for some $\beta\in\Br(Y)$, then $Y$ is a supersingular K3 surface.
\end{Proposition}
\begin{proof}
  The derived category of any twisted Fourier-Mukai partner of a twisted K3 surface must have trivial Serre functor. Furthermore, the Betti numbers of a surface are twisted derived invariants. Over the complex numbers, this follows for instance from \cite[Remark 5.30]{Huy06}. Identical reasoning gives a proof in any characteristic, upon replacing singular cohomology with $l$-adic or crystalline cohomology. It follows from the classification of surfaces that any twisted Fourier-Mukai partner $(Y,\beta)$ of the twisted K3 surface $(X,\alpha)$ is again a twisted K3 surface. As the isocrystal $\tH(X/K)$ has slope 1, the same is true for $\tH(Y/K)$. It follows that $Y$ is supersingular.
\end{proof}
By analogy with the complex case, one expects that $\Phi^{\cry}_{v^{-B\boxplus B'}(P^{\bullet})}$ in fact preserves the integral structure. This is the content of the following result.
\begin{Theorem}\label{thm:Mainprelim}
  Let $(X,\alpha)$ and $(Y,\beta)$ be twisted supersingular K3 surfaces with B-field lifts $B$ and $B'$. If $P^\bullet\in D(X\times Y,-\alpha\boxtimes\beta)$ is a complex of twisted sheaves such that the induced Fourier-Mukai functor
  \[
    \Phi_{P^\bullet}\colon D(X,\alpha)\to D(Y,\beta)
  \]
  is an equivalence, then the induced cohomological transform gives an isomorphism
    \[
      \Phi^{\cry}_{v^{-B\boxplus B'}(P^{\bullet})}:(\tH(X/W,B),\tN(X,B))\iso(\tH(Y/W,B'),\tN(Y,B'))
    \]
    of marked K3 crystals.
\end{Theorem}
As it is somewhat unrelated to the rest of this paper, we defer the proof to Appendix \ref{app:twistedChernCharacters} and Appendix \ref{app:deformation}. An immediate consequence of Theorem \ref{thm:Mainprelim} is that the Artin invariant is preserved under a derived equivalence.
\begin{Corollary}\label{cor:FMpartner2}
  If $(X,\alpha)$ and $(Y,\beta)$ are twisted supersingular K3 surfaces that are derived equivalent, then $\sigma_0(X,\alpha)=\sigma_0(Y,\beta)$.
\end{Corollary}

We mention a few consequences for supersingular K3 surfaces with small and large Artin invariants.
\begin{Corollary}
  The unique supersingular K3 surface of Artin invariant 1 has no non-trivial twisted Fourier-Mukai partners.
\end{Corollary}
\begin{proof}
Any twisted supersingular K3 surface $(X,\alpha)$ with $\alpha\neq 0$ must have Artin invariant $\sigma_0\geq 2$.  So, every twisted partner is in fact untwisted, and has Artin invariant 1.
\end{proof}

\begin{Corollary}
  If $X$ is a supersingular K3 surface of Artin invariant 10 and $\alpha\in\Br(X)$ is a nontrivial Brauer class, then $(X,\alpha)$ is not derived equivalent to any non-twisted variety.
\end{Corollary}
\begin{proof}
  By Proposition \ref{prop:FMpartner1}, any potential Fourier-Mukai partner must be a supersingular K3 surface. Because $\alpha$ is nontrivial, $\sigma_0(X,\alpha)=11$. But, any non-twisted supersingular K3 surface has Artin invariant $\sigma_0\leq 10$.
\end{proof}

Over the complex numbers, it is known that the action of a derived equivalence on the singular Mukai lattice preserves a certain orientation structure (see \cite{MR2553878}, and \cite{2017arXiv171100846R} for the twisted case). As we do not have access to the singular cohomology lattice in positive characteristic, we instead formulate the orientation as a structure on the extended N\'{e}ron-Severi group.

Let $N$ be a lattice of signature $(s_+,s_-)$.
\begin{Definition}\label{def:orientation}
  An \textit{orientation} on $N$ is a choice of orientation (in the sense of vector spaces) on a positive definite subspace $P\subset N_{\bR}$ of dimension $s_+$. An isometry $g\colon N\to N'$ or $N_{\bQ}\to N'_{\bQ}$ of oriented lattices is \textit{orientation preserving} if the composition
  \[
    P\to N_{\bR}\xrightarrow{g_{\bR}}N'_{\bR}\to P'
  \]
  is orientation preserving in the usual sense, where $N'_{\bR}\to P'$ is the orthogonal projection.
\end{Definition}

We refer to \cite{2017arXiv171100846R} for further discussion. Let $X$ be a K3 surface. The extended N\'{e}ron-Severi group $\tN(X)$ has signature $(2,22)$. Let $H\in\Pic(X)$ be an ample class. We consider the positive definite subspace $P_H=\langle (1,0,-1),H\rangle\subset\tN(X)\otimes\bR$. By the given ordering of the basis elements, $P_H$ is endowed with an orientation in the sense of vector spaces, and hence $\tN(X)$ acquires an orientation in the sense of Definition \ref{def:orientation}. Note that this orientation does not depend on the choice of $H$, in the sense that the oriented lattices resulting from different choices are isometric by an orientation preserving isometry (namely the identity).

Suppose that $(X,\alpha)$ is a twisted K3 surface and $\alpha$ is killed by $p$. Given any B-field lift $B$ of $\alpha$, we obtain a natural orientation on $\tN(X,B)$ by declaring that the isometry
\[
  e^B\colon \tN(X,B)_{\bQ}\iso\tN(X)_{\bQ}
\]
is orientation preserving. If $P^{\bullet}\in D(X\times Y,-\alpha\boxtimes\beta)$ is a complex of twisted sheaves that induces a Fourier-Mukai equivalence $\Phi_{P^{\bullet}}\colon D(X,\alpha)\to D(Y,\beta)$, then by Theorem \ref{thm:Mainprelim} the cohomological transform $\Phi^{\cry}_{v^{-B\boxplus B'}(P^{\bullet})}$ restricts to an isometry $\tN(X,B)\iso\tN(Y,B')$. We say that $\Phi^{\cry}_{v^{-B\boxplus B'}(P^{\bullet})}$ is \textit{orientation preserving} if this isometry is orientation preserving in the sense of Definition \ref{def:orientation}.\footnote{We do not really need Theorem \ref{thm:Mainprelim} to make this definition, as it is immediate that the cohomological transform gives an isometry after tensoring with $\bQ$, which is enough to define the property of being orientation preserving.}
\begin{Conjecture}\label{conj:mainPrelim2}
  With the assumptions of Theorem \ref{thm:Mainprelim}, the cohomological transform $\Phi^{\cry}_{v^{-B\boxplus B'}(P^{\bullet})}$ is orientation preserving (in the sense that the induced map on twisted N\'{e}ron-Severi groups is orientation preserving).
\end{Conjecture}
We resolve this conjecture in some special cases in Appendix \ref{app:deformation}.

\subsection{Crystalline Torelli theorems for twisted supersingular K3 surfaces}

In this section we formulate a crystalline Torelli theorem relating isomorphisms between the twisted Mukai crystals associated to a twisted supersingular K3 surface to Fourier-Mukai equivalences between derived categories of twisted sheaves. We will first prove a result relating the twisted Mukai crystal to equivalences of abelian categories of twisted sheaves. Suppose that $(X,\alpha)$ and $(Y,\beta)$ are twisted K3 surfaces. If $f\colon X\iso Y$ is an isomorphism such that $f^*\beta=\alpha$, then for any B-field lifts $B,B'$ of $\alpha$ and $\beta$, the induced map $f_*\colon\tH(X/K)\iso\tH(Y/K)$ restricts to an isomorphism
\[
  \tH(X/W,B)\iso\tH(Y/W,B')
\]
of K3 crystals. We wish to identify which such maps of crystals are induced by an isomorphism of surfaces. Consider the natural filtration on $\tH(X/K)$ by codimension. If $B$ is any B-field lift of a Brauer class $\alpha$, then there is an induced filtration $0\subset F^2\subset F^1\subset F^0=\tH(X/W,B)$ on the submodule $\tH(X/W,B)\subset\tH(X/K)$ defined by
\begin{align*}
  F^2&=\tH(X/W,B)\cap\H^4(X/K)\\
  F^1&=\tH(X/W,B)\cap(\H^2(X/K)\oplus\H^4(X/K))\\
  F^0&=\tH(X/W,B)\cap(\H^0(X/K)\oplus\H^2(X/K)\oplus\H^4(X/K))
\end{align*}
Note that the $F^i$ are in fact subcrystals. There is a canonical isomorphism
\[
  F^1/F^2\iso\H^2(X/W)
\]
induced by $(0,h,h.B)\mapsto h$. This isomorphism is compatible with the Frobenius operators on both sides. Furthermore, as $F^2$ is isotropic, the Mukai pairing induces a pairing on the quotient $F^1/F^2$, and under the above isomorphism this pairing is taken to the cup product pairing.

\begin{Definition}
  We say that an isomorphism $g\colon\tH(X/W,B)\to\tH(Y/W,B')$ of K3 crystals \textit{preserves the codimension filtrations} if $g(F^i)\subset F^i$ for $i=0,1,2$.
\end{Definition}
In fact, because $F^1$ is exactly the orthogonal complement of $F^2$, this is equivalent to $g(F^2)\subset F^2$. If $g$ preserves the filtrations, then it induces a map
\[
  \frac{F^1\tH(X/W,B)}{F^2\tH(X/W,B)}\iso\frac{F^1\tH(Y/W,B')}{F^2\tH(Y/W,B')}
\]
and hence an isomorphism $g_0\colon\H^2(X/W)\iso\H^2(Y/W)$ of K3 crystals.
\begin{Theorem}\label{thm:twistedOgus1}
  Let $(X,\alpha)$ and $(Y,\beta)$ be twisted supersingular K3 surfaces. Let $B$ and $B'$ be B-field lifts of $\alpha$ and $\beta$. If
  \[
    g\colon (\tH(X/W,B),\tN(X,B))\to(\tH(Y/W,B'),\tN(Y,B'))
  \]
  is an isomorphism of marked K3 crystals such that
  \begin{enumerate}
      \item $g$ sends $(0,0,1)$ to $(0,0,1)$ (and in particular preserves the codimension filtrations), and
      \item $g_0$ maps an ample class to an ample class,
  \end{enumerate}
  then there exists an isomorphism $f\colon X\to Y$ satisfying $f^*\beta=\alpha$ and a line bundle $\sL$ on $Y$ such that $\exp(c_1^{\cry}(\sL)/p)\circ f_*=g$.
\end{Theorem}
\begin{proof}
  As we have observed above, the assumption that $g$ is a map of K3 crystals implies that the induced isomorphism $g_0\colon\H^2(X/W)\iso\H^2(Y/W)$ is also a map of K3 crystals. Furthermore, as $g$ restricts to an isomorphism between the respective N\'{e}ron-Severi groups, $g_0$ restricts to an isomorphism between the respective (non-twisted) N\'{e}ron-Severi groups, and by assumption sends an ample class to an ample class. Thus, by Ogus's crystalline Torelli theorem (in the formulation of Theorem \ref{thm:Ogus1}), $g_0$ is induced by a unique isomorphism $f\colon X\iso Y$.
  
  Consider the commuting diagram
  \[
    \begin{tikzcd}
      \tH(X/W,B)\arrow{rr}{g}\arrow{dr}[swap]{f_*}&&\tH(Y/W,B')\\
      &\tH(Y/W,g_0(B))\arrow{ur}[swap]{\exp(B'-g_0(B))}&
    \end{tikzcd}
  \]
  of $W$-modules. As both $g$ and $f_*$ are isomorphisms of marked K3 crystals, so is $\exp(B'-g_0(B))$. The inclusion $\tN(Y,g_0(B))\subset\tN(Y)\otimes\bQ$ becomes an isomorphism after tensoring with $\bQ$, so $\exp(B'-g_0(B))\in\tN(Y)\otimes\bQ$. It follows that $B'-g_0(B)=c_1^{\cry}(\sL)/p$ for some line bundle $\sL$ on $Y$. This completes the proof.

\end{proof}


We next prove the main result of this paper: a Torelli theorem relating isomorphisms of the twisted Mukai crystals to equivalences of derived categories. This should be seen as a supersingular analog of the twisted derived Torelli theorem of Huybrechts and Stellari \cite{MR2310257}.
\begin{Theorem}\label{thm:Main1}
  Let $(X,\alpha)$ and $(Y,\beta)$ be twisted supersingular K3 surfaces with B-field lifts $B$ and $B'$.
  \begin{enumerate}
    \item If $\Phi_{P^\bullet}:D(X,\alpha)\to D(Y,\beta)$ is a Fourier-Mukai equivalence, the induced map
    \[
      \Phi^{\cry}_{v^{-B\boxplus B'}(P^\bullet)}:(\tH(X/W,B),\tN(X,B))\iso(\tH(Y/W,B'),\tN(Y,B'))
    \]
    is an isomorphism of marked K3 crystals (which if Conjecture \ref{conj:mainPrelim2} holds is orientation preserving).
    \item If
    \[
      g:(\tH(X/W,B),\tN(X,B))\iso(\tH(Y/W,B'),\tN(Y,B'))
    \]
    is an orientation preserving isomorphism of marked K3 crystals, then there exists a Fourier-Mukai equivalence $\Phi_{P^\bullet}:D(X,\alpha)\to D(Y,\beta)$ such that $\Phi^{\cry}_{v^{-B\boxplus B'}(P^\bullet)}=g$.
  \end{enumerate}
\end{Theorem}
\begin{proof}
Part (1) is Theorem \ref{thm:Mainprelim}. Let us show (2). The reader will note that the proof is essentially the same as the proof of the main theorem of \cite{MR2310257}.  We start with an isomorphism
\[
  g:(\tH(X/W,B),\tN(X,B))\iso(\tH(Y/W,B'),\tN(Y,B'))
\]
of oriented marked crystals. By definition, this means that we have a diagram
\[
  \begin{tikzcd}
    \tH(X/W,B)\arrow{r}{g}&\tH(Y/W,B')\\
    \tN(X,B)\arrow[hook]{u}\arrow{r}{g}&\tN(Y,B')\arrow[hook]{u}
  \end{tikzcd}
\]
where the upper horizontal arrow is an isomorphism of K3 crystals and the vertical arrows are the canonical inclusions. The lower horizontal arrow is then an isometry with respect to the Mukai pairings on the twisted N\'{e}ron-Severi lattices, and preserves the canonical orientations. Let $w=(r,l,s)$ be the image $g(0,0,1)$ of the class of a point under $g$. We will first treat the case when $r\neq 0$. We may then assume that $r>0$ by composing with $g_0=-\id_{\tH(Y/W,B')}$, which lifts to the shift functor $F\mapsto F\left[1\right]$. As $g$ is an isometry we have $g^2=0$, so the moduli space $\sM_{(Y,\beta)}(w)$ of $\beta$-twisted sheaves on $Y$ with twisted Mukai vector $w$ is a $\bG_m$-gerbe over a K3 surface $M_{(Y,\beta)}(w)$, and there is a $B$-field $B''\in H^2(M_{(Y,\beta)}(w)/W)$ and an equivalence $D(Y,\beta)\cong D(M_{(Y,\beta)}(w),\alpha_{B''})$ inducing an orientation preserving isometry $\tH(Y/W,B')\to\tH(M_{(Y,\beta)}(w)/W,B'')$ (the relevant facts about moduli spaces of twisted sheaves have been extended to our setting in Section 4 of \cite{BL17}). The inverse of this equivalence sends $(0,0,1)$ to $w$, so composing we may assume without loss of generality that $g(0,0,1)=(0,0,1)$.

To apply Theorem \ref{thm:twistedOgus1}, we need to ensure that the induced map $g_0\colon\H^2(X/W)\to\H^2(Y/W)$ maps an ample class to an ample class (see the discussion preceding the statement of Theorem \ref{thm:twistedOgus1}). By results of \cite{Ogus83} there exists a sequence of $(-2)$ classes $C_1,\dots,C_n$ such that composing $g_0$ with the composition of reflections $C_n\circ\dots\circ C_1$ takes the ample cone $C_X$ of $X$ to $\pm C_Y$. But $g$ was assumed to be orientation preserving, so in fact $C_X\mapsto C_Y$. As explained in \cite{MR2310257}, each $C_i$ is induced by a spherical twist functor $T_{C_i}\colon D(Y,\beta)\to D(Y,\beta)$ on the derived category, and moreover the action of such a functor on the twisted Mukai crystal sends $(0,0,1)$ to $(0,0,1)$. We have reduced to the case when the assumptions of Theorem \ref{thm:twistedOgus1} apply. We therefore find an isomorphism $f\colon X\to Y$ satisfying $f^*\beta=\alpha$ and a line bundle $\sL$ on $Y$ such that
\[
  \exp(c_1^{\cry}(\sL)/p)\circ f_*=g
\]
Both of these maps lift to Fourier-Mukai equivalences, so this completes the proof.

\end{proof}

If we are only interested in the existence of a derived equivalence, it turns out that we may forget the marking and orientation. This will make our computations more manageable.

\begin{Proposition}\label{prop:forgetmarking}
   Let $(H,N)$ and $(H',N')$ be marked oriented supersingular K3 crystals. If there exists an isomorphism $H\iso H'$ of K3 crystals (not necessarily preserving the marking or orientation), then there exists an isomorphism $(H,N)\iso (H',N')$ of marked K3 crystals that is orientation preserving.
\end{Proposition}
\begin{proof}
   By results of Nikulin (see for instance Theorem 15.1.5 of \cite{Huy06}), there exists an isometry $N\iso N'$. The remainder of the proof follows that of \cite[Theorem 7.4]{Ogus78} exactly.
\end{proof}

\begin{Theorem}\label{thm:main2}
  Let $(X,\alpha)$ and $(Y,\beta)$ be twisted supersingular K3 surfaces, with B-field lifts $B,B'$ of $\alpha$ and $\beta$. The following are equivalent.
  \begin{enumerate}
      \item There exists a derived equivalence $D(X,\alpha)\cong D(Y,\beta)$.
      \item There exists an isomorphism $\tH(X/W,B)\cong\tH(Y/W,B')$ of K3 crystals (not necessarily preserving any marking or orientation).
  \end{enumerate}
\end{Theorem}
\begin{proof}
  This follows from Theorem \ref{thm:Main1} combined with Proposition \ref{prop:forgetmarking}.
\end{proof}

\section{Applications to Fourier-Mukai equivalences}

\subsection{Characteristic subspaces}\label{sec:characteristicsubspaces}
It is a remarkable discovery of Ogus that supersingular K3 crystals are determined by certain semi-linear algebraic data, called characteristic subspaces. These are very accessible for making explicit computations. This correspondence resembles the equivalence between Hodge structures and periods. In this section we will explain this connection, and reformulate the theorems of the previous section in terms of characteristic subspaces.

Let $V$ be a vector space over $\bF_p$ of even dimension $2\sigma_0$ equipped with a non-degenerate, non-neutral bilinear form. Let $\sigma\colon k\to k$ be the Frobenius map $\lambda\mapsto\lambda^p$, and define 
\[
  \varphi=\id\otimes\,\sigma\colon V\otimes k\to V\otimes k
\]
\begin{Definition}
  A \textit{characteristic subspace} of $V$ is a subspace $K\subset V\otimes k$ such that
  \begin{enumerate}
    \item $K$ is totally isotropic of dimension $\sigma_0$ and
    \item $K+\varphi(K)$ has dimension $\sigma_0+1$.
  \end{enumerate}
  If in addition $\sum_{i=0}^\infty\varphi^i(K)=V\otimes k$, we say that $K$ is \textit{strictly characteristic}.
\end{Definition}

\begin{Definition}
  A \textit{characteristic subspace datum} is a pair $(K,V)$, where $V$ is a vector space over $\bF_p$ of even dimension equipped with a non-degenerate, non-neutral bilinear form, and $K\subset V\otimes k$ is a strictly characteristic subspace.
  
  An isomorphism $(K,V)\iso (K',V')$ of characteristic subspace data is an isometry $g:V\to V'$ such that $g\otimes k$ maps $K$ to $K'$.
\end{Definition}

The \textit{Artin invariant} of a characteristic subspace datum is the dimension of $K$ as a $k$-vector space. Let $H$ be a supersingular K3 crystal with Tate module $T$ and Artin invariant $\sigma_0$. The quotient $T_0=T^\vee/T$ is an $\bF_p$-vector space of dimension $2\sigma_0$, which inherits a non-degenerate and non-neutral quadratic form from $T$, defined by
\[
  \overline{v}.\overline{w}=p^{-1}\langle pv,pw\rangle
\]
for $v,w\in T^{\vee}$. There are natural inclusions
\[
  T\otimes W\subset H\subset T^{\vee}\otimes W
\]
By Theorem 3.20 of \cite{Ogus78}, the image $K$ of $H$ in $(T^{\vee}\otimes W)/(T\otimes W)\cong T_0\otimes k$ is a strictly characteristic subspace.
\begin{Definition}\label{def:definition of the char sub of a crystal}
  The characteristic subspace datum associated to a supersingular K3 crystal $H$ is the pair $(K,T_0)$.
\end{Definition}
\begin{Remark}
  It is perhaps more correct to replace $K$ with $\varphi^{-1}(K)$, which is also strictly characteristic (this is the convention taken in \cite{Ogus78}). This distinction is important when working over non-perfect bases.
\end{Remark}

Any isomorphism $H\iso H'$ of supersingular K3 crystals restricts to an isometry $T\iso T'$ of Tate modules, and hence induces an isomorphism $(K,T_0)\iso (K',T'_0)$ of characteristic subspace data. By results of Ogus, the converse is also true.

\begin{Proposition}\label{prop:char subspaces and crystals}
  Suppose that $H,H'$ are supersingular K3 crystals of rank $n$. Let $T\subset H$ and $T'\subset H'$ be their Tate modules, and $(K,T_0)$ and $(K',T'_0)$ the corresponding characteristic subspace data. Any isomorphism $(K,T_0)\iso (K',T'_0)$ of characteristic subspace data is induced by an isomorphism $H\iso H'$ of K3 crystals.
\end{Proposition}
\begin{proof}
  Consider an isomorphism $(K,T_0)\to (K',T'_0)$ of characteristic subspace data. By Corollary 3.18 of \cite{Ogus78}, there exists an isometry $T_0\iso T'_0$. Moreover, by results of Nikulin (see for instance Theorem 14.2.4 of \cite{Huy06}), every isometry $T_0\iso T'_0$ lifts to an isometry $T\iso T'$. The result then follows from Theorem 3.20 of \cite{Ogus78}.
\end{proof}

\begin{Notation}\label{not:subspacedatum}
  Let $X$ be a supersingular K3 surface. We write $T(X)$ for the Tate module of $\H^2(X/W)$ and write $K(X)$ for its associated characteristic subspace. Similarly, we write $\tT(X)$ and $\tK(X)$ for the Tate module and characteristic subspace associated to $\tH(X/W)$.
  
  If $\alpha\in\Br(X)$ is a Brauer class and $B$ is a B-field lift of $\alpha$, then $\tH(X/W,B)$ is a supersingular K3 crystal of rank 24, and we let $\tT(X,B)$ be its Tate module and $\tK(X,B)$ its associated characteristic subspace.
\end{Notation}
\begin{Remark}
  Suppose that $X$ is a supersingular K3 surface. We can give an alternate description of the characteristic subspace datum associated to $\H^2(X/W)$ that does not use crystalline cohomology. Consider the first Chern class map
  \begin{equation}\label{eq:first chern class map}
    c_1^{dR}\colon\Pic(X)\to\H^2_{dR}(X/k)
  \end{equation}
  Tensoring with $k$, we get a map
  \begin{equation}\label{eq:chern}
    N(X)\otimes k\to\H^2_{dR}(X/k)
  \end{equation}
  As the discriminant group of $N(X)$ is killed by $p$, there is a natural inclusion $pN(X)^{\vee}\subset N(X)$. This induces an inclusion
  \[
    \dfrac{pN(X)^{\vee}}{pN(X)}\subset\dfrac{N(X)}{pN(X)}
  \]
  of $\bF_p$-vector spaces. Multiplication by $p$ induces an isomorphism
  \[
    N(X)_0\defeq\dfrac{N(X)^{\vee}}{N(X)}\iso\dfrac{pN(X)^{\vee}}{pN(X)}
  \]
  This gives a natural inclusion
  \[
    N(X)_0\otimes k\hto N(X)\otimes k
  \]
  Unwinding definitions, we see that the composition
  \[
    K(X)\hto T(X)_0\otimes k=N(X)_0\otimes k\hto N(X)\otimes k 
  \]
  identifies $K(X)$ with the kernel of (\ref{eq:chern}).
\end{Remark}
We now translate the crystalline Torelli theorems in previous sections into the language of characteristic subspaces. Let $X$ be a supersingular K3 surface. We write $(K(X),T(X)_0)$ for the characteristic subspace datum associated to the supersingular K3 crystal $\H^2(X/W)$ (see Notation \ref{not:subspacedatum}). Using Proposition \ref{prop:char subspaces and crystals}, Ogus's crystalline Torelli theorem \ref{thm:Ogus2} can be rephrased in terms of characteristic subspaces.
\begin{Theorem}\label{thm:characteristictorelli1}
  If $X$ and $Y$ are supersingular K3 surfaces, then $X$ is isomorphic to $Y$ if and only if there is an isomorphism $(K(X),T(X)_0)\iso (K(Y),T(Y)_0)$ of characteristic subspace data.
\end{Theorem}
Suppose that $\alpha\in\Br(X)$ is a Brauer class with B-field lift $B$. The supersingular K3 crystal $\tH(X/W,B)$ gives rise to the characteristic subspace datum $(\tK(X,B),\tT(X,B)_0)$. If $\alpha$ is non-trivial, then the functional on $\tT(X,B)$ given by pairing with $(0,0,1)$ is divisible by $p$. Hence, $p^{-1}(0,0,1)$ gives rise to a non-trivial element of $\tT(X,B)_0$. Anticipating future notation, let us write $f_X$ for the element $p^{-1}(0,0,1)\in\tT(X,B)_0$. The twisted crystalline Torelli theorem then assumes the following form.
\begin{Theorem}\label{thm:characteristictorelli2}
  Let $X$ and $Y$ be supersingular K3 surfaces and $\alpha\in\Br(X),\beta\in\Br(Y)$ non-trivial Brauer classes with B-field lifts $B,B'$. There exists an isomorphism $(X,\alpha)\cong (Y,\beta)$ if and only if there exists an isomorphism of characteristic subspace data $(\tK(X,B),\tT(X,B)_0)\iso (\tK(Y,B'),\tT(Y,B')_0)$ sending $f_X$ to $f_Y$.
\end{Theorem}

Finally, we restate our main theorem on derived equivalences in terms of characteristic subspaces.
\begin{Theorem}\label{thm:characteristictorelli3}
  Let $X$ and $Y$ be supersingular K3 surfaces and $\alpha\in\Br(X),\beta\in\Br(Y)$ Brauer classes with B-field lifts $B,B'$. There exists a Fourier-Mukai equivalence $D(X,\alpha)\cong D(Y,\beta)$ if and only if there exists an isomorphism $(\tK(X,B),\tT(X,B)_0)\iso(\tK(Y,B'),\tT(Y,B')_0)$ of characteristic subspace data.
\end{Theorem}
\begin{proof}
  This follows immediately from Proposition \ref{prop:char subspaces and crystals} and Theorem \ref{thm:main2}.
\end{proof}

We summarize our results in the following table. In each row, the object in the left column is classified by the indicated K3 crystal, and by the indicated characteristic subspace datum. We continue the notation $f_X=p^{-1}(0,0,1)$. 

\vskip 1\baselineskip

\newcommand\Tstrut{\rule{0pt}{2.6ex}}         
\newcommand\Bstrut{\rule[-0.9ex]{0pt}{0pt}}   

\begin{center}
\begin{tabular}{ c|c|c } 
  \bf{Geometric Object} & \bf{K3 Crystal} & \bf{Characteristic Subspace Datum}\Bstrut\\
  \hline
  $X$           & $H^2(X/W)$                    & $(K(X),T(X)_0)$\Tstrut \\ 
  $(X,\alpha)$  & $\tH(X/W,B)$ + $(0,0,1)\in\tH$& $(\tK(X,B),\tT(X,B)_0)$ + $f_X\in \tT(X,B)_0$   \\
  $D(X)$        & $\tH(X/W)$                    & $(\tK(X),\tT(X)_0)$    \\
  $D(X,\alpha)$ & $\tH(X/W,B)$                  & $(\tK(X,B),\tT(X,B)_0)$
\end{tabular}
\end{center}

\vskip 1\baselineskip

Let $X$ be a supersingular K3 surface and $B$ a B-field on $X$. We will examine the relationship between the characteristic subspace data associated to $\H^2(X/W)$ and $\tH(X/W,B)$. This relationship is particularly simple if we make a special choice of B-field lift. We may characterize Artin's subgroup $\U^2(X,\m_p)\subset\H^2(X,\m_p)$ as follows.
\begin{Definition}\label{def:alternatedef}
  The subgroup $\U^2(X,\m_p)\subset\H^2(X,\m_p)$ consists of exactly those classes in $\H^2(X,\m_p)$ which, under the inclusion $\H^2(X,\m_p)\hto\H^2_{dR}(X/k)$, are orthogonal to the image of the first Chern class map~\eqref{eq:first chern class map}.
\end{Definition}
\begin{Definition}
  A B-field $B=\frac{a}{p}\in\H^2(X/K)$ is \textit{transcendental} if the image of $a$ in $\H^2_{dR}(X/k)$ is contained in $\U^2(X,\m_p)$.
\end{Definition}
By the above characterization of $\U^2(X,\m_p)$, we see that a B-field $B=\frac{a}{p}$ is transcendental if and only if $\langle a,D\rangle$ is divisible by $p$ for every $D\in N(X)$, or equivalently if and only if $B\in T(X)^{\vee}\otimes W\subset p^{-1}\H^2(X/W)$.

Now, it also follows from the above that the image of $\U^2(X,\m_p)\hto\H^2_{dR}(X/k)$ is contained in the subspace $V\otimes k/K\subset\H^2_{dR}(X/k)$, where $V=N(X)_0$. We may identify the image of this map as follows.
\begin{Lemma}\label{lem:coolbijection}
  The inclusion $\U^2(X,\m_p)\hto\H^2_{dR}(X/k)$ induces an isomorphism 
  \[
    \U^2(X,\mu_p)\cong\dfrac{\left\{B\in V\otimes k|B-\varphi(B)\in K+\varphi(K)\right\}}{K}\subset \dfrac{V\otimes k}{K}
  \]
\end{Lemma}
\begin{proof}
  This is proved in \cite{BL17}. See Remark 3.3.9ff and Lemma 3.3.15.
\end{proof}
Equivalently, this map may be described in the following way: given a class $\alpha\in\U^2(X,\m_p)$, we choose a B-field $B=\frac{a}{p}$ such that $\overline{a}=d\log\alpha$. Such a $B$ is transcendental, so it is contained in $T(X)^{\vee}\otimes W$. We send $\alpha$ to the image of $B$ under the composition
\[
  T(X)^{\vee}\otimes W\to\dfrac{T(X)^\vee\otimes W}{T(X)\otimes W}=T(X)_0\otimes k=V\otimes k\to \dfrac{V\otimes k}{K}
\]
Under the isomorphism of Lemma \ref{lem:coolbijection}, the subgroup $pN(X)^{\vee}/pN(X)$ described by (\ref{eq:artinssubgroup}) is identified with the fixed points of $\varphi$ acting on $V\otimes k$.

Suppose $\alpha\in\Br(X)$ is a non-trivial Brauer class. Let $B=\frac{a}{p}$ be a transcendental B-field lift of $\alpha$. One can show that $\tT(X,B)$ is then the free $W$-module
\[
  \tT(X,B)=\langle (p,0,0),(0,D,0),(0,0,1)\rangle
\]
where $D$ ranges over all classes in $N(X)$ (see Proposition 3.4.21 of \cite{BL17}). There is an induced decomposition
\[
  \tT(X,B)_0=T(X)_0\oplus U_2 
\]
where $U_2$ is the copy of the hyperbolic plane generated by the vectors $e$ and $f$, given respectively by pairing with $p^{-1}(p,0,0)$ and $p^{-1}(0,0,1)$. If $x_1,\dots,x_{\sigma_0}$ is a basis for $K(X)$, then the characteristic subspace $\tK(X,B)$ associated to $\tH(X/W,B)$ by Definition \ref{def:definition of the char sub of a crystal} is given by
\begin{equation}\label{eq:decomposition222}
  \tK(X,B)=\left\langle x_1+(x_1.B)f,\dots,x_{\sigma_0}+(x_{\sigma_0}.B)f,e+B+\dfrac{B^2}{2}f\right\rangle
\end{equation}

With the above descriptions as motivation, we record the following result on abstract characteristic subspaces. We fix an $\bF_p$-vector space $V$ of dimension $2\sigma_0$, equipped with a non-degenerate, non-neutral bilinear form, and a characteristic subspace $K\subset V\otimes k$. Let $U_2$ be a copy of the hyperbolic plane over $\bF_p$, equipped with the standard basis $v,w$ satisfying $v^2=w^2=0$ and $v.w=-1$. Set $\tV=V\oplus U_2$. This is a vector space of dimension $2\sigma_0+2$, whose natural bilinear form is again non-degenerate and non-neutral.

\begin{Proposition}\label{prop:radbijection}
Let $x_1,\dots,x_{\sigma_0}$ be a basis for $K$. The map
 \[
  \dfrac{\left\{B\in V\otimes k|B-\varphi(B)\in K+\varphi(K)\right\}}{K}\to\left\{
  \begin{tabular}{@{}c@{}}
  characteristic subspaces $\tK\subset\tV\otimes k$ \\
  such that $v\notin\tK$ and $\tK\cap v^\perp/v=K$\\
  \end{tabular}
  \right\}
  \]
  given by
  \[
    B\mapsto \left\langle x_1+(x_1.B)v,\dots,x_{\sigma_0}+(x_{\sigma_0}.B)v,w+B+\frac{B^2}{2}v\right\rangle
  \]
  is a bijection.
 \end{Proposition}
 \begin{proof}
   This is prove in Proposition 3.1.11 of \cite{BL17} (see also Remark 3.2.3).
 \end{proof}

\subsection{Fourier-Mukai partners of twisted supersingular K3 surfaces}\label{sec:applications}
We will now deduce some consequences for derived equivalences of twisted supersingular K3 surfaces. We have attached to a supersingular K3 surface two supersingular K3 crystals: the second crystalline cohomology $\H^2(X/W)$, which has rank 22, and the Mukai crystal $\tH(X/W)$, which has rank 24. Note that there is an inclusion $\H^2(X/W)\subset\tH(X/W)$, which induces an equality $T(X)_0=\tT(X)_0$ identifying the corresponding characteristic subspaces. Thus, the characteristic subspace data of the two crystals are \textit{equal} (of course, the crystals have different ranks, and thus are not isomorphic). Combining this observation with Theorem \ref{thm:characteristictorelli1} and Theorem \ref{thm:characteristictorelli3}, we obtain another proof of the following result of Lieblich and Olsson \cite{LO15}.

\begin{Theorem}\label{thm:Untwisted1}
  If $X$ and $Y$ are supersingular K3 surfaces, then $D(X)\cong D(Y)$ if and only if $X\cong Y$.
\end{Theorem}
\begin{proof}
  If $D(X)\cong D(Y)$, then by Theorem \ref{thm:Main1} we have in particular an isomorphism $\tH(X/W)\cong\tH(Y/W)$ of K3 crystals. Thus, there is an isomorphism $(\tK(X),\tT(X)_0)\cong(\tK(Y),\tT(Y)_0)$ of characteristic subspace data. But by the above observation, this gives an isomorphism $(K(X),T(X)_0)\cong (K(Y),T(Y)_0)$ between the characteristic subspace data associated to $\H^2(X/W)$ and $\H^2(Y/W)$. By Theorem \ref{thm:characteristictorelli2}, this implies $X\cong Y$.
  
  The converse is immediate.
\end{proof}

\begin{Remark}
  The characteristic subspace associated to a K3 crystal might be viewed as an analog of the transcendental lattice. Theorem \ref{thm:Untwisted1} is then analogus to the classical result that two complex K3 surfaces with Picard number $\geq 12$ are derived equivalent if and only if they are isomorphic (see Corollary 16.3.7 of \cite{Huy06}).
\end{Remark}

\begin{Corollary}\label{cor:Untwisted2}
  A twisted supersingular K3 surface $(X,\alpha)$ is derived equivalent to a non-twisted K3 surface if and only if $\sigma_0(X,\alpha)\leq 10$.
\end{Corollary}
\begin{proof}
  A non-twisted supersingular K3 surface has $\sigma_0\leq 10$. So, if $\sigma_0(X,\alpha)=11$ then $(X,\alpha)$ cannot be derived equivalent to a non-twisted K3 surface.
  
  Conversely, suppose that $\sigma_0(X,\alpha)\leq 10$. Consider the characteristic subspace datum $(\tK(X,B),\tT(X,B)_0)$ associated to $\tH(X/W,B)$. By the surjectivity of the period morphism, we find a supersingular K3 surface $Y$ and an isomorphism $(K(Y),T(Y)_0)\cong (\tK(X,B),\tT(X,B)_0)$ of characteristic subspace data. But $(K(Y),T(Y)_0)=(\tK(Y),\tT(Y)_0)$, so Theorem \ref{thm:characteristictorelli3} implies that $D(Y)\cong D(X,\alpha)$.
\end{proof}

\begin{Proposition}
  Let $E$ be a supersingular elliptic curve, and let $X=\Km(E\times E)$ be the unique supersingular K3 surface of Artin invariant $\sigma_0=1$. Every supersingular K3 surface of Artin invariant 2 is derived equivalent to $D(X,\alpha)$ for some Brauer class $\alpha\in\Br(X)$.
\end{Proposition}
\begin{proof}
  Let $Y$ be a supersingular K3 surface of Artin invariant 2. Let $(K(Y),T(Y)_0)$ be the characteristic subspace datum associated to $\H^2(Y/W)$. Pick a non-zero isotropic vector $v\in T(Y)_0$ (such a vector exists because the dimension of $T(Y)_0$ is $\geq 4$). Set $V=v^\perp$. We may then find an isotropic vector $w\in T(Y)_0$ such that $v.w=-1$, and such that there is an orthogonal decomposition $T(Y)_0=V\oplus U_2$ where $U_2$ is the hyperbolic plane generated by $v,w$. As $K(Y)$ is strictly characteristic, $v$ is not in $K(Y)$. Let $K=K(Y)\cap v^{\perp}/v$. This is a strictly characteristic subspace of $V$. By the surjectivity of the period morphism, the pair $(K,V)$ is isomorphic to the characteristic subspace datum associated to $\H^2(X/W)$, where $X$ is the unique supersingular K3 surface of Artin invariant $1$. By Proposition \ref{prop:radbijection} and Lemma \ref{lem:coolbijection}, there is a Brauer class $\alpha\in\Br(X)$ and a transcendental B-field lift $B$ of $\alpha$ such that $(K(X,B),T(X,B)_0)\cong(K(Y),T(Y)_0)$. As before, $(K(Y),T(Y)_0)$ is equal to the characteristic subspace datum associated to $\tH(Y/W)$, and Theorem \ref{thm:characteristictorelli3} then gives the result.
\end{proof}

We also prove the following (compare to Corollary 7.9 of \cite{HS04}).
\begin{Proposition}\label{thm:caldararusthing}
  If $X$ is a supersingular K3 surface and $\alpha\in\Br(X)$ is a Brauer class, then for any integer $k$ such that $(k,p)=1$, there is a derived equivalence
  \[
    D(X,\alpha)\cong D(X,\alpha^k)
  \]
\end{Proposition}
\begin{proof}
  Suppose that $\alpha$ is non-trivial. Let $B$ be a transcendental B-field lift of $\alpha$. We have a canonical decomposition
  \[
    \tT(X,B)_0=T(X)_0\oplus U_2
  \]
  where $U_2$ is a copy of the hyperbolic plane generated by $e,f$, corresponding to pairing with $p^{-1}(p,0,0)$ and $p^{-1}(0,0,1)$ respectively. If $x_1,\dots,x_{\sigma_0}$ is a basis for $K(X)$, then we have that
  \[
    \tK(X,B)=\left\langle x_1+(x_1.B)f,\dots,x_{\sigma_0}+(x_{\sigma_0}.B)f,e+B+\dfrac{B^2}{2}f\right\rangle
  \]
  (see the notation and discussion preceding (\ref{eq:decomposition222})). Let $\lambda\in\bF_p^{\times}$ be the image of $k$ modulo $p$. Consider the isometry $m_{\lambda}\colon\tT(X,B)_0\iso\tT(X,B)_0$ defined by $v\mapsto v$ for $v\in T(X)_0$, $e\mapsto \lambda^{-1}e$, and $f\mapsto \lambda f$. The image of $\tK(X,B)$ under $m_{\lambda}$ is
  \[
    \left\langle x_1+(x_1.\lambda B)f,\dots,x_{\sigma_0}+(x_{\sigma_0}.\lambda B)f,e+\lambda B+\dfrac{(\lambda B)^2}{2}f\right\rangle
  \]
  The element $\lambda B$ is again in $\U^2(X,\m_p)$, and its image in the Brauer group is $\alpha^k$ (here, we are abusing notation by identifying $B$ with its image in $\U^2(X,\m_p)\subset\H^2_{dR}(X/k)$, and with its image under the isomorphism of Lemma \ref{lem:coolbijection}ff). If $\lambda'\in\bZ_p^{\times}$ is a lift of $\lambda$, then $\lambda'B$ is a (transcendental) B-field lift of $\alpha^k$, and $m_{\lambda}$ gives an isomorphism $(\tK(X,B),\tT(X,B)_0)\iso (\tK(X,\lambda' B),\tT(X,\lambda' B)_0)$ of characteristic subspace data. We conclude the result by Theorem \ref{thm:characteristictorelli3}.
\end{proof}


As a further application of our main theorem, we will count the number of twisted Fourier-Mukai partners of a twisted supersingular K3 surface. This will require some more involved calculations with characteristic subspaces.


\subsection{Interlude: the orthogonal group of a strictly characteristic subspace}\label{sec:orthogonal}

Let $V$ be a vector space over $\bF_p$ of dimension $2\sigma_0$, equipped with a non-degenerate, non-neutral bilinear form, and let $K\subset V\otimes k$ be a strictly characteristic subspace.
\begin{Definition}
   We define $O_{K}(V)$ to be the automorphism group of the characteristic subspace datum $(K,V)$. That is, $O_K(V)$ is the group of isometries $g:V\to V$ such that $g\otimes k$ maps $K$ to itself.
\end{Definition}
In this section we will study the orthogonal group $O_K(V)$ using certain coordinates introduced by Ogus. Because $K$ is strictly characteristic, the subspace
\[
  l_K=K\cap\varphi(K)\cap\dots\cap\varphi^{\sigma_0-1}(K)
\]
is one-dimensional, and
\[
  \varphi^{\sigma_0-1}(K)=l_K+\varphi(l_K)+\dots+\varphi^{\sigma_0-1}(l_K)
\]
Pick a generator $e$ for $L$, and scale $e$ so that $e.\varphi^{\sigma_0}(e)=1$. Such a choice of $e$ is unique up to multiplication by a $p^{\sigma_0}+1$-th root of unity. Set $e_i=\varphi^{i-1}(e)$, so that $\left\{e_1,\dots,e_{\sigma_0}\right\}$ is a basis for $\varphi^{\sigma_0-1}(K)$ and $\left\{e_1,\dots,e_{2\sigma_0}\right\}$
is a basis for $V\otimes k$. In this basis both the bilinear form and the Frobenius take a particularly simple form. Define
\begin{equation}\label{eq:structureconstants}
  a_i=e_1.e_{\sigma_0+i+1}\quad 1\leq i\leq\sigma_0-1
\end{equation}
By Lemma 3.22 of \cite{Ogus78}, the bilinear form on $V\otimes k$ is given by the matrix
\begin{equation}\label{eqn:form}
  \begin{pmatrix} 
    0 & A \\
    A^T & 0 
  \end{pmatrix}
\end{equation}
where $A$ is the $\sigma_0\times\sigma_0$-matrix
\[
  \begin{pmatrix} 
    1      & a_1    &         a_2  &         a_3    & \dots  & a_{\sigma_0-1}   \\
    0      & 1      &       F(a_1) &       F(a_2)   & \dots  & F(a_{\sigma_0-2})  \\
    0      & 0      &         1    &       F^2(a_1) & \dots  & F^2(a_{\sigma_0-3})\\
    0      & 0      &         0    & 1              & \dots  & F^3(a_{\sigma_0-4})\\
    \vdots & \vdots & \vdots       & \vdots         & \ddots & \vdots                   \\
    0      & 0      & 0            & 0              & \dots  & 1
  \end{pmatrix}
\]
The Frobenius is given by $\varphi(e_i)=e_{i+1}$ for $1\leq i\leq 2\sigma_0-1$, and
\[
  \varphi(e_{2\sigma_0})=\lambda_1e_1+\dots+\lambda_{\sigma_0}e_{\sigma_0}+\mu_1e_{\sigma_0+1}+\dots+\mu_{\sigma_0}e_{2\sigma_0}
\]
for some particular scalars $\lambda_i,\mu_i$ which satisfy $\lambda_1=1$ and $\mu_1=0$; the rest can be described in terms of the $a_i$. In this basis, every element of $O_K(V)$ is diagonalized. Indeed, suppose that $g\in O_K(V)$ is an isometry that preserves $K$. Because $g$ commutes with $\varphi$, it must also preserve $l_K$, and therefore
\[
  g(e)=\zeta e
\]
for some scalar $\zeta$. Because of the uniqueness of $e$, $\zeta$ must be a $(p^{\sigma_0}+1)$-th root of unity. We have $g(e_i)=\zeta^{p^{i-1}}(e_1)$, so $g$ is the diagonal matrix
\begin{equation}\label{eqn:diagonal}
  \begin{pmatrix} 
    \zeta &           &               &        & \\
            & \zeta^p &               &        & \\
            &           & \zeta^{p^2} &        & \\
            &           &               & \ddots & \\
            &           &               &        &\zeta^{p^{2\sigma_0-1}}\\
  \end{pmatrix}
\end{equation}
In particular, $g$ is determined by $\zeta$. The map $g\mapsto\zeta$ defines an injective homomorphism
\begin{equation}\label{eqn:inclusion}
  O_K(V)\hto \mu_{p^{\sigma_0}+1}(k)
\end{equation}

We will determine the image of (\ref{eqn:inclusion}). 

\begin{Proposition}\label{prop:subgroup}
  Let $m$ be a divisor of $\sigma_0$ such that $\sigma_0/m$ is odd. If $\zeta$ is a primitive $p^m+1$-th root of unity, then $\zeta$ is contained in the image of (\ref{eqn:inclusion}) if and only if $a_i=0$ for every $1\leq i\leq\sigma_0-1$ such that $2m$ does not divide $i$.
\end{Proposition}
\begin{proof}
  Let $m$ be a divisor of $\sigma_0$ with $\sigma_0/m$ odd, let $\zeta$ be a primitive $(p^m+1)$-th root of unity, and let $g$ be the corresponding diagonal element (\ref{eqn:diagonal}) of $\GL(V\otimes k)$.
  
  Suppose that $\zeta$ is contained in the image of (\ref{eqn:inclusion}). As $g$ is an isometry, we have for each $1\leq i\leq\sigma_0-1$ that
  \[
    a_i=e_1.e_{\sigma_0+i+1}=g(e_1).g(e_{\sigma_0+i+1})=\zeta^{p^{\sigma_0+i}+1}a_i
  \]
  If $a_i\neq 0$ for some $1\leq i\leq\sigma_0-1$, then $\zeta^{p^{\sigma_0+i}+1}=1$, so $p^m+1$ divides $p^{\sigma_0+i}+1$. Thus, $p^{\sigma_0+i}\equiv -1$ modulo $p^m+1$. The order of $p$ in $(\bZ/(p^m+1)\bZ)^\times$ is $2m$, and because $\sigma_0/m$ is odd, we have $p^{\sigma_0}\equiv -1$ modulo $p^m+1$. It follows that $p^i\equiv 1$ modulo $p^m+1$, and therefore $2m$ divides $i$.
  
  Conversely, suppose that $a_i=0$ for every $1\leq i\leq\sigma_0-1$ such that $2m$ does not divide $i$. We will show that $g\in O_K(V)$. It is clear that any matrix of the form (\ref{eqn:diagonal}) preserves $K$. We next check that $g\in O(V\otimes k)$, that is, that $g$ preserves the bilinear form given by the matrix (\ref{eqn:form}). We will show that
  \[
    g(e_k).g(e_j)=e_k.e_j
  \]
  for all $1\leq k,j\leq 2\sigma_0$. It is clear that $\varphi(g(e_i))=g(\varphi(e_i))$ for $1\leq i\leq 2\sigma_0-1$, so it will suffice to check that
  \begin{equation}\label{eq:misceq2}
     g(e_1).g(e_{\sigma_0+i+1})=e_1.e_{\sigma_0+i+1}=a_i
  \end{equation}
  for all $1\leq i\leq\sigma_0-1$. We have
  \[
     g(e_1).g(e_{\sigma_0+i+1})=\zeta^{p^{\sigma_0+i}+1}a_i
  \]
  If $2m$ does not divide $i$, then by assumption $a_i=0$, so (\ref{eq:misceq2}) holds. Suppose that $2m$ divides $i$. It follows that $p^i\equiv 1$ modulo $p^m+1$. We know that $p^{\sigma_0}\equiv -1$ modulo $p^m+1$, so $p^{\sigma_0+i}+1\equiv 0$ modulo $p^m+1$. We conclude that $\zeta^{p^{\sigma_0+i}+1}=1$. It follows that $g$ is an isometry.
  
  Finally, we will show that $g$ is defined over $\bF_p$, that is, that it is in the image of the map $O(V)\to O(V\otimes k)$. Equivalently, we will show that $g$ commutes with the Frobenius $\varphi$. As previously observed, $\varphi(g(e_i))=g(\varphi(e_i))$ for $1\leq i\leq 2\sigma_0-1$, so it remains only to verify that $\varphi(g(e_{2\sigma_0}))=g(\varphi(e_{2\sigma_0}))$. The pairing is non-degenerate, so it will suffice to show that
  \begin{equation}\label{eq:misceq1}
    \varphi(g(e_{2\sigma_0})).e_k=g(\varphi(e_{2\sigma_0})).e_k
  \end{equation}
  for all $2\leq k\leq 2\sigma_0+1$ (the vectors $e_2,\dots,e_{2\sigma_0+1}$ are a basis for $V\otimes k$). Using that $g$ is an isometry, (\ref{eq:misceq1}) becomes
  \[
    \zeta^{p^{2\sigma_0}+p^{k-1}}\varphi(e_{2\sigma_0}).e_k=\varphi(e_{2\sigma_0}).e_k
  \]
  We compute that
  \[
    \varphi(e_{2\sigma_0}).e_k=
    \begin{cases}
    F^{k-1}(a_{\sigma_0-k+1})&\mbox{ if }2\leq k\leq\sigma_0\\
    1&\mbox{ if }k=\sigma_0+1\\
    0&\mbox{ if }\sigma_0+2\leq k\leq 2\sigma_0+1
    \end{cases}
  \]
  It follows that (\ref{eq:misceq1}) is an equality for $\sigma_0+1\leq k\leq 2\sigma_0+1$. Suppose that $2\leq k\leq\sigma_0$. If $2m$ does not divide $\sigma_0-k+1$, then by assumption $a_{\sigma_0-k+1}=0$, and so (\ref{eq:misceq1}) is an equality. Suppose that $2m$ divides $\sigma_0-k+1$. Because $m$ divides $\sigma_0$, this implies that $m$ divides $2\sigma_0-k+1$. Furthermore, as $\sigma_0/m$ is odd, so is $(2\sigma_0-k+1)/m$. It follows that $p^m+1$ divides $p^{2\sigma_0-k+1}+1$, and therefore $p^m+1$ divides $p^{2\sigma_0}+p^{k-1}$. Thus, (\ref{eq:misceq1}) is an equality in this case as well. We conclude that $g\in O_K(V)$, and hence $\zeta$ is in the image of (\ref{eqn:inclusion}).

  
\end{proof}

\begin{Proposition}
  The image of (\ref{eqn:inclusion}) is equal to the subgroup $\mu_{p^m+1}(k)\subset\mu_{p^{\sigma_0}+1}(k)$, where $m$ is the largest divisor of $\sigma_0$ such that
  \begin{enumerate}
      \item $\sigma_0/m$ is odd, and
      \item if $1\leq i\leq\sigma_0-1$ and $2m$ does not divide $i$ then $a_i=0$.
  \end{enumerate}
  or $m=0$ if no such divisor exists.
\end{Proposition}
\begin{proof}
  The inclusion (\ref{eqn:inclusion}) realizes $O_K(V)$ as a subgroup of the cyclic group $\mu_{p^{\sigma_0}+1}(k)$. Thus, $O_K(V)$ is itself cyclic, say of order $k$. Let $g$ be a generator, and let $\zeta$ be the image of $g$ in $\mu_{p^{\sigma_0}+1}(k)$. Let $m'$ be the smallest non-negative integer such that $k$ divides $p^{m'}+1$. The order of $p$ modulo $k$ is then $2m'$. We also know that $p^{\sigma_0}+1$ is divisible by $k$, so therefore $m'$ divides $\sigma_0$, and $\sigma_0/m'$ is odd. By Proposition \ref{prop:subgroup} applied to $m'$, we deduce that $a_i=0$ for every $1\leq i\leq\sigma_0-1$ such that $2m'$ does not divide $i$. That is, $m'$ is a divisor of $\sigma_0$ satisfying properties (1) and (2). Furthermore, as $g$ was a generator, we have that the image of (\ref{eqn:inclusion}) is contained in $\mu_{p^{m'}+1}(k)$. On the other hand, it follows from Proposition \ref{prop:subgroup} that $\mu_{p^m+1}(k)$ is contained in the image of (\ref{eqn:inclusion}). We conclude that $\mu_{p^m+1}(k)\subset\mu_{p^{m'}+1}(k)$. By maximality of $m$, we deduce that $m=m'$, and the result follows.
\end{proof}
Note that if the structure constants $a_i$ associated to $K$ are all non-zero then $m=0$. Thus, if $K$ is generic then $\#O_K(V)=2$. We remark that a similar calculation is carried out in Theorem 3.1 of \cite{2017arXiv170202290J}. There is also given a table of the dimensions of the locus of characteristic subspaces with different values of $m$ for each Artin invariant $1\leq \sigma_0\leq 10$.

\subsection{Counting twisted Fourier-Mukai partners}\label{sec:counting}


We apply our techniques to count twisted partners of supersingular K3 surfaces.  The analogous computation over the complex numbers is carried out in \cite{MR2550163}.  We fix a twisted supersingular K3 surface $(X,\alpha)$ with Artin invariant $\sigma_0=\sigma_0(X,\alpha)$. We define
\[
  \FM(X,\alpha)=\left\{
  \begin{tabular}{@{}c@{}}
  twisted supersingular K3 surfaces $(Y,\beta)$ \\
  such that $D(X,\alpha)\cong D(Y,\beta)$\\
  \end{tabular}
\right\}/\,\mbox{isomorphism}
\]
Note that by Proposition \ref{prop:FMpartner1}, any twisted partner of $(X,\alpha)$ must be a twisted supersingular K3 surface. Our goal is to enumerate this set. We will do this by using characteristic subspaces and the computations of Section \ref{sec:orthogonal}.

\begin{Definition}
  If $V$ is an $\bF_p$-vector space with a bilinear form, we let $I(V)$ denote the set of non-zero vectors $v\in V$ such that $v^2=0$.
\end{Definition}
There is an obvious action of $O_K(V)$ on $I(V)$. Fix a B-field lift $B$ of $\alpha$. Let $(K,V)=(\tK(X,B),\tT(X,B)_0)$ be the characteristic subspace datum associated to $\tH(X/W,B)$.

\begin{Theorem}\label{thm:Count1}
  \begin{equation}\label{eq:count}
    \#\FM(X,\alpha)=
      \begin{cases}
        \# O_{K}(V)\setminus I(V)+1&\mbox{ if }\sigma_0\leq 10\\
        \# O_{K}(V)\setminus I(V)&\mbox{ if }\sigma_0=11\\
      \end{cases}
  \end{equation}
\end{Theorem}
We let $\FM^d(X,\alpha)$ be the subset of $\FM(X,\alpha)$ consisting of those Fourier-Mukai partners $(X_1,\alpha_1)$ such that $\alpha_1$ has order $d$.  Of course, in our situation $\FM^d(X,\alpha)$ is empty unless $d=1,p$.  We can dispense with the $d=1$ case quickly.
\begin{Lemma}
  If $\sigma_0\leq 10$, then $\#\FM^1(X,\alpha)=1$. If $\sigma_0=11$, then $\#\FM^1(X,\alpha)=0$.
\end{Lemma}
\begin{proof}
  This follows immediately from Corollary \ref{cor:Untwisted2}.
\end{proof}

We next turn to $\FM^p(X,\alpha)$. Note that if $\sigma_0=1$, then this set is empty, so let us assume that $\sigma_0\geq 2$. If $(Y,\beta)\in \FM^p(X,\alpha)$, then there exists a $B$ field lift $B'$ of $\beta$ and an isometry 
\[
  g:\tT(Y,B')_0\iso V=\tT(X,B)_0
\]
that sends $\tK(Y,B')$ to $K=\tK(X,B)$. Because the order of $\beta$ is $p$, the functional given by pairing with the vector $(0,0,1)\in\tT(Y,B')$ is divisible by $p$. As before, we denote by $f_Y$ the resulting (non-trivial) element $p^{-1}(0,0,1)\in\tT(Y,B')_0$. The image of this element under $g$ gives an element of $I(V)$. Note that this element does not depend on our choice of $B'$, although it does depend on $g$. We define a map of sets
\[
  \mu:\FM^p(X,\alpha)\to O_K(V)\setminus I(V)
\]
by $\mu(\left[(Y,\beta)\right])=\left[g(f_Y)\right]$, where $O_K(V)$ acts on $I(V)$ on the left in the obvious way. This is well defined, independent of the choices of $B'$ and $g$.
\begin{Proposition}
  The map $\mu$ is a bijection.
\end{Proposition}
\begin{proof}
  We first show injectivity. Consider two Fourier-Mukai partners 
  \[
    (X_1,\alpha_1),(X_2,\alpha_2)\in\FM^p(X,\alpha)
  \]
  such that $\mu(X_1,\alpha_1)=\mu(X_2,\alpha_2)$. Pick $B$-field lifts $B_1,B_2$, and choose isometries $g_i:\tT(X_i,B_i)_0\iso\tT(X,B)_0$, $i=1,2$ that preserve the respective characteristic subspaces. By assumption, there exists $g\in O_K(V)$ such that
  \[
    g(g_1(f_{X_1}))=g_2(f_{X_2})
  \]
  Thus, the isometry
  \[
    g_2^{-1}\circ g\circ g_1:\tT(X_1,B_1)_0\to\tT(X_2,B_2)_0
  \]
  preserves the characteristic subspaces and satisfies $f_{X_1}\mapsto f_{X_2}$. By Theorem \ref{thm:characteristictorelli2}, we find an isomorphism $(X_1,\alpha_1)\cong(X_2,\alpha_2)$.
  
  We next show surjectivity. Pick a non-zero isotropic vector $v\in V$. The vector space $V^-=v^\perp/v$ has dimension $2\sigma_0-2$, and because we assumed $\sigma_0\geq 2$, its induced bilinear form is again non-degenerate and non-neutral. Consider the subspace 
  \[
    K^-=(K\cap v^\perp)/v\subset V^-\otimes k
  \]
  Pick an isotropic vector $w\in V$ such that $v.w=-1$ and a decomposition
  \[
    V=V^-\oplus\langle v,w\rangle
  \]
  Let $\left\{x_1,\dots,x_{\sigma_0-1}\right\}$ be a basis for $K^-$.  By Proposition \ref{prop:radbijection}, there is an element $B'\in V^-\otimes k$ such that
  \[
    K=\left\langle x_1+(x_1.B')v,\dots,x_{\sigma_0-1}+(x_{\sigma_0-1}.B')v,w+B'+\frac{B'^2}{2}w\right\rangle
  \]
  By the surjectivity of the period morphism, the pair $(K^-,V^-)$ is isomorphic to the characteristic subspace datum associated to $\H^2(Y/W)$ for some supersingular K3 surface $Y$. Under this isomorphism, $B'$ corresponds via Lemma \ref{lem:coolbijection} to some class in $\U^2(Y,\m_p)$. Let $\beta$ be its image in $\Br(Y)$. By abuse of notation, let $B'$ denote a B-field lift of $\beta$ that maps to the given class in $\U^2(Y,\m_p)$. By construction, we have an isomorphism $(\tK(Y,B'),\tT(Y,B')_0)\iso (K,V)$ sending $f_Y$ to $v$. By Theorem \ref{thm:characteristictorelli3}, $(Y,\beta)$ is a Fourier-Mukai partner of $(X,\alpha)$, and its image under $\mu$ is $v$, as desired.
\end{proof}

Using the results of Section \ref{sec:orthogonal}, we can explicitly determine the size of the set
\[
  O_K(V)\setminus I(V)
\]
in terms of the structure constants $a_i$ associated to $K$.
\begin{Lemma}\label{lem:free}
  The action of $O_K(V)$ on $I(V)$ is free. That is, each $v\in I(V)$ is fixed only by the identity.
\end{Lemma}
\begin{proof}
  If $g\in O_K(V)$ fixes $v$, then
  \[
    v.e=g(v).g(e)=v.(\zeta e)=\zeta(v.e)
  \]
  Because $v$ is non-zero and fixed by $\varphi$, we must have $v.e\neq 0$.  It follows that $\zeta=1$, and hence $g$ is the identity.
\end{proof}

  Let $(X,\alpha)$ be a twisted supersingular K3 surface with Artin invariant $\sigma_0$, and let $a_i$ be be the constants associated to the strictly characteristic subspace $\tK(X,\alpha)$ corresponding to $(X,\alpha)$ (see \ref{eq:structureconstants}). Let $m$ be the largest divisor of $\sigma_0$ such that $\sigma_0/m$ is odd and $a_i=0$ for all $1\leq i\leq\sigma_0-1$ such that $2m$ does not divide $i$, or set $m=0$ if no such divisor exists.
\begin{Theorem}\label{Count2}
  The number of twisted Fourier-Mukai partners of $(X,\alpha)$ is given by the formula
  \[
    \#\FM(X,\alpha)=
      \begin{dcases}
        \frac{p^{\sigma_0}+1}{p^m+1}(p^{\sigma_0-1}-1)+1 &\mbox{ if }1\leq\sigma_0\leq 10\\
        \frac{p^{\sigma_0}+1}{p^m+1}(p^{\sigma_0-1}-1)   &\mbox{ if }\sigma_0=11\\
      \end{dcases}
  \]
\end{Theorem}
\begin{proof}
  By eg. \cite[Lemma 4.12ff]{Ogus78}, if $V$ is a vector space over $\bF_p$ of dimension $2\sigma_0$ equipped with a non-degenerate, non-neutral bilinear form, then
  \[
    \# I(V)=(p^{\sigma_0}+1)(p^{\sigma_0-1}-1)
  \]
  By Lemma \ref{lem:free}, the action of $O_K(V)$ on $I(V)$ is free. The result follows from Theorem \ref{thm:Count1}.
\end{proof}
\begin{Remark}
  The integer $m$ associated to $(X,\alpha)$ is equal to 0 if all of the structure constants $a_i$ are non-zero. This is the case for instance if $X$ is a generic supersingular K3 surface, where by generic we mean that the corresponding point in Ogus's moduli space of marked supersingular K3 surfaces \cite{Ogus83} is contained in the complement of an explicit proper closed subset.
\end{Remark}





\appendix

\section{Twisted Chern characters}\label{app:twistedChernCharacters}
In this section we will discuss the twisted Chern characters introduced in Definitions \ref{def:twistedcherncharacter} and \ref{def:twisted Chern character part 2}. 

\subsection{Cocycle formulation}\label{ssec:cocycles}

We begin by reformulating our definition in terms of cocycles. Let $X$ be a smooth projective variety over an algebraically closed field of arbitrary characteristic and $\alpha\in\Br(X)$ a Brauer class. Choose a sufficiently fine \'{e}tale cover $U_i\to X$ and a cocycle $\alpha_{ijk}\in \mathbb{G}_m(U_{ijk})$ representing $\alpha$. An $\alpha$-twisted sheaf $\ms E$ is a tuple $(\ms E_i,\varphi_{ij})$ where each $\ms E_i$ is a coherent sheaf on $U_i$ and the $\varphi_{ij}$ are isomorphisms on the intersections satisfying $\varphi_{ij}\circ\varphi_{jk}\circ\varphi_{ki}=\alpha_{ijk}$. As $n\alpha=0$, the cocycle $\alpha^n_{ijk}$ represents the zero class in cohomology. However, identifying the category of twisted sheaves defined with respect to the cocycle $\alpha^n_{ijk}$ with the category of sheaves on $X$ requires a choice. This choice is provided by our lift $\alpha'$ of $\alpha$. To see this, consider a cocycle $\alpha'_{ijk}\in \mu_n(U_{ijk})$ representing $\alpha'$. Using the inclusion $\mu_n\subset\mathbb{G}_m$, we may view this also as a cocycle for $\alpha$. This cocycle has the property that $\alpha'^n_{ijk}=1$. Let us use this cocycle to represent the category of $\alpha$-twisted sheaves. The tuple $(\ms E_i^{\otimes n},\varphi_{ij}^{\otimes n})$ satisfies the cocycle condition, and hence glues to a sheaf on $X$. Moreover, we have
\[
    \ch^{\alpha'}(\ms E)=\sqrt[n]{\ch(\ms E_i^{\otimes n},\varphi_{ij}^{\otimes n})}
\]

\subsection{Comparison with the twisted Chern character of Huybrechts and Stellari}\label{ssec:comparison}

Suppose that the ground field is the complex numbers. Huybrechts and Stellari introduce in \cite[Section 1]{HS04} a twisted Chern character for twisted sheaves on $X$, whose definition relies on the existence of an invertible twisted sheaf in the differentiable category. Thus, their definition is a priori not algebraic. Nevertheless, we will show that it is essentially equivalent to our notion, up to some differences in convention. We remark that a comparison between the twisted Chern characters of \cite{HS04} and others in the literature is given in \cite{MR2310257}.

The Chern character of a locally free sheaf is defined by using the splitting principle to reduce to the case of line bundles.
In this paper, as well as \cite{BL17,LMS11}, the convention used is to define the first Chern character modulo $n$ of a line bundle $L$ on $X$ to be equal to the image of the class $[L]\in\H^1(X,\mathbb{G}_m)$ under the boundary map
\[
    \delta:\H^1(X,\mathbb{G}_m)\to\H^2(X,\mu_n)
\]
induced by the Kummer sequence. In \cite[Proposition 1.2]{HS04}, Huybrechts and Stellari (implicitly) use the convention under which the first Chern character of $L$ is $c_1(L)=-\delta(L)$ (there are similar sign conventions according to the various targets of the different incarnations of the first Chern class map). This difference in conventions will be the cause of a sign discrepancy in our comparison. In what follows, we will (as in the rest of this paper) use the normalization $c_1(L)=\delta(L)$. Using $\ch^{\vee}$ to denote the Chern character taken according to the normalization of \cite{HS04}, we have
\[
    \ch^{\vee}(\mathscr{E})=\ch(\mathscr{E}^{\vee})
\]

Let $B$ be a B-field lift of $\alpha$ in the sense of \cite{HS04}. Thus, $B$ is a class in $\H^2(X,\mathbb{Q})$ whose image under the exponential map is $\alpha$. Let $B_{ijk}\in\mathbb{Q}(U_{ijk})$ be a cocycle representative for $B$, and use the cocycle $\alpha_{ijk}=\exp(B_{ijk})$ for represent $\alpha$. Let $\mathscr{E}=(\mathscr{E}_i,\varphi_{ij})$ be a locally free $\alpha$--twisted sheaf. To avoid confusion, write $\ch_{\HS}^B(\mathscr{E})$ for the twisted Chern character of \cite[Section 1]{HS04} evaluated on $\mathscr{E}$ (denoted by $\ch^B(\mathscr{E})$ in loc. cit.). This is defined by choosing real valued differentiable functions $a_{ij}$ such that $B_{ijk}=-a_{ij}-a_{jk}-a_{ki}$. The tuple $(\mathscr{E}_i,\varphi_{ij}\exp(a_{ij}))$ then satisfies the cocyle condition, and hence glues to a sheaf $\mathscr{E}_B$ on $X$, and Huybrechts and Stellari define
\[
    \ch_{\HS}^B(\mathscr{E})=\ch(\mathscr{E}_B^{\vee})
\]
(here, as described above, the right hand side is the Chern character taken with the normalization used in the rest of this paper). Let $n$ be a positive integer such that $nB\in\H^2(X,\mathbb{Z})$. Consider the class $\alpha'\in\H^2(X,\mu_n)$ given by the image of $nB$ under the composition
\[
    \H^2(X,\mathbb{Z})\twoheadrightarrow\H^2(X,\mathbb{Z}/n\mathbb{Z})\xrightarrow[\exp]{\sim}\H^2(X,\mu_n)
\]
The twisted Chern characters of Definition \ref{def:twisted Chern character part 2} are expressed in terms of those of \cite{HS04} by
\[
    \ch^{\alpha'}(\mathscr{E})=\exp(-B)\ch^{-B}_{\HS}(\mathscr{E}^{\vee})
\]
This follows from the cocycle formulation in~\ref{ssec:cocycles} by taking $n$th powers of both sides.

\subsection{Integrality of twisted Chern characters}\label{ssec:integrality}

We now prove Proposition \ref{prop:integrality1} and Theorem \ref{thm:Mainprelim}. Although we have phrased the results in this section for $p$-torsion Brauer classes and crystalline B-fields, it will be clear that the proofs apply essentially unchanged to the singular cohomology constructions of \cite{HS04}, the $l$-adic \'{e}tale theory of \cite{LMS11}, and the de Rham-Witt theory for $p^n$-torsion Brauer classes.

We will begin by comparing twisted sheaves on $\bG_m$-gerbes, $\m_n$-gerbes, and Brauer-Severi varieties. Let $X$ be a smooth projective variety over $k$. We assume that the characteristic of $k$ is a prime $p$, possibly equal to $2$. Let $\pi\colon\sX\to X$ be a $\m_n$-gerbe. If $\sE$ is a locally free coherent sheaf on $\sX$, we define
  \[
    \ch_{\sX}(\sE)=\sqrt[\leftroot{-2}\uproot{2}n]{\ch(\pi_*(\sE^{\otimes n}))}\in A^*(X)\otimes\bQ
  \]
We compare this twisted Chern character to the notion introduced in Definition \ref{def:twistedcherncharacter}. Let $\sX_{\bG_m}$ be the $\bG_m$-gerbe obtained from $\sX$. The inclusion $\m_n\hto\bG_m$ induces a natural map $\iota\colon \sX\to\sX_{\bG_m}$, and $\iota_*$ induces an equivalence between the respective categories of twisted sheaves. The gerbe $\sX$ gives a $\m_n$-gerbe lift of $\sX_{\bG_m}$, and hence gives rise to a choice of invertible $n$-fold twisted sheaf $\sL$ on $\sX_{\bG_m}$ (see Definition \ref{def:twistedcherncharacter}ff). This sheaf satisfies $\iota^*\sL\cong\cO_{\sX}$.
It follows that for any locally free twisted sheaf $\sE$ on $\sX$, we have
\[
  \ch_{\sX}(\sE)=\ch^{\sL}(\iota_*\sE)
\]
The pushforward $\sA=\pi_*\sEnd(\sE)$ is an Azumaya algebra on $X$ whose associated cohomology class in $\H^1(X,\PGL_n)$ maps to the Brauer class of $\sX$ under the boundary map $\H^1(X,\PGL_n)\to\H^2(X,\bG_m)$. Let $f\colon P\to X$ be the Brauer-Severi scheme associated to $\sA$. We have a Cartesian diagram
\[
  \begin{tikzcd}
    \bP_{\sX}(\sE)\arrow{r}{\pi'}\arrow{d}[swap]{g}&P\arrow{d}{f}\\
    \sX\arrow{r}{\pi}&X
  \end{tikzcd}
\]
The projective bundle $\bP_{\sX}(\sE)\to\sX$ comes with a universal quotient $g^*\sE\twoheadrightarrow\cO(1)$. We may also view $\bP_{\sX}(\sE)$ as a $\m_n$-gerbe via $\pi'$, and the invertible sheaf $\cO(1)$ is 1-twisted with respect to this gerbe structure. Define an invertible sheaf $\cO(n)=\pi'_*\cO(1)^{\otimes n}$ on $P$. The $\m_n$-gerbe $\bP_{\sX}(\sE)\to P$ is essentially trivial (that is, has trivial Brauer class), and is isomorphic to the gerbe of $n$-th roots of $\cO(n)$. Under this isomorphism, the universal quotient $\cO(1)$ acquires a second universal property: it is the universal (twisted) invertible sheaf equipped with an isomorphism $\cO(1)^{\otimes n}\iso \pi'^*\cO(n)$.



\begin{Proposition}\label{prop:cohomology}
  There is a natural isomorphism of graded $W$-algebras $\H^*(P/W)\cong\H^*(X/W)[x]/f(x)$ where $x\in\H^2(P/W)$ is a class restricting to $c_1^{\cry}(\cO(1))$ on geometric fibers, and $f(x)$ is a monic polynomial of degree $n$. The map $f^*\colon \H^*(X/W)\otimes k\to\H^*(P/W)\otimes k$ is injective.
\end{Proposition}
\begin{proof}
  This follows from the Leray spectral sequence associated to $P\to X$ (see for instance \cite[Lemma 1.6]{Yos06}).
\end{proof}
\begin{Proposition}\label{prop:integral1}
    Suppose that $X$ is a K3 surface, $\alpha\in\Br(X)$ is a Brauer class of order $p$, and $\sE$ is an $\alpha$-twisted locally free sheaf. If $B$ is a B-field lift of $\alpha$, then $\ch^{B}(\sE)\in\tN(X,B)$.
\end{Proposition}
\begin{proof}
  Let $\pi\colon\sX\to X$ be a $\m_p$-gerbe whose associated cohomology class maps to $\alpha$. We regard $\sE$ as a twisted sheaf on $\sX$. It will suffice to show the result for a particular choice of B-field lift. Consider a B-field lift $B=\frac{a}{p}$ of $\alpha$ where the image of $a$ modulo $p$ is equal to $d\log$ of the cohomology class $[\sX]\in\H^2(X,\m_p)$. It follows that $f^*B=c_1^{\cry}(\cO(p))/p+h$ for some $h\in\H^2(P/W)$.
  
  By the above discussion, we have $\ch_{\sX}(\sE)=\ch^B(\iota_*\sE)$. The Chern character of a coherent sheaf on a K3 surface is integral, so $\ch_{\sX}(\sE)^p\in\tN(X)$. Thus, $\ch_{\sX}(\sE)\in\tN(X)\otimes\bZ[p^{-1}]$. To prove the result, it will therefore suffice to show that $e^{-B}\ch_{\sX}(\sE)\in\H^*(X/W)$. By Proposition \ref{prop:cohomology}, this is true if and only if $f^*(e^{-B}\ch_{\sX}(\sE))\in\H^*(P/W)$. We compute
  \begin{align*}
    f^*(e^{-B}\ch_{\sX}(\sE))&=f^*(e^{-B}).f^*\ch_{\sX}(\sE)\\
    &=e^{-h}.\sqrt[\leftroot{-2}\uproot{2}p]{\ch(-c_1^{\cry}(\cO(p)).\ch(f^*\pi_*(\sE^{\otimes p}))}\\
    &=e^{-h}.\sqrt[\leftroot{-2}\uproot{2}p]{\ch(-c_1^{\cry}(\cO(p)).\ch(\pi'_*g^*(\sE^{\otimes p}))}\\
    &=e^{-h}.\sqrt[\leftroot{-2}\uproot{2}p]{\ch(\pi'_*(\pi'^*(\cO(-p))\otimes g^*(\sE^{\otimes p})))}\\
    &=e^{-h}.\sqrt[\leftroot{-2}\uproot{2}p]{\ch(\pi'_*((g^*(\sE)\otimes\cO(-1))^{\otimes p})))}\\
    &=e^{-h}.\ch(\pi'_*(g^*(\sE)\otimes\cO(-1)))
  \end{align*}
Write $G=\pi'_*(g^*(\sE)\otimes\cO(-1))$. 
We have
\[
  e^{-h}\ch(G)=(r,c_1(G)-rh,c_1(G)^2/2+c_2(G)-h.c_1(G)+rh^2/2,0,\dots,0)
\]
where $r=\rk(G)$. The Chern classes $c_i(G)$ are all integral. In particular, this already gives the result if $p\neq 2$. Suppose $p=2$. We compute 
\[
  c_1(G)=c_1(\pi_*\det(\sE))-f^*(a)+ph
\]
where $a=pB$. It is well known that the pairing on the N\'{e}ron-Severi group of a K3 surface is even. By \cite[Theorem 4.7]{Ogus83}, the pairing on $\H^2(X/W)$ is also even, in the sense that $2$ divides $x.x$ for any $x\in\H^2(X/W)$ (this is a non-trivial statement only if $p=2$). We obtain that $c_1(G)^2$ is divisible by 2, as is $r$, which gives the result.
\end{proof}
Proposition \ref{prop:integrality1} follows, as the Todd class of a K3 surface is integral. Let $(X,\alpha)$, $(Y,\beta)$ be twisted K3 surfaces such that $\alpha$ and $\beta$ are killed by $p$, and fix B-field lifts $B,B'$. Define
\begin{equation}\label{eq:MoLatticesMoProblems}
  \tN(X\times Y,-B\boxplus B')=(N^*(X\times Y)\otimes\bZ[p^{-1}])\cap(e^{-B\boxplus B'}H^*(X\times Y))
\end{equation}
Using the method of Proposition \ref{prop:integral1}, we can show that the twisted Chern characters of a twisted sheaf on $X\times Y$ lie in this lattice. Unfortunately, we must assume here that $p\geq 5$, although we would be surprised if this assumption was necessary. 
In any case, under the additional assumptions that $X$ and $Y$ are supersingular and that $\Phi_{\sE}$ is an equivalence we will prove the result for all $p$ in Appendix \ref{app:deformation}.
\begin{Proposition}\label{prop:integral2}
If $\sE\in D(X\times Y,-\alpha\boxtimes\beta)$ is a locally free twisted sheaf and $p\geq 5$, then
  \[
    \ch^{-B\boxplus B'}(\sE)\in\tN(X\times Y,-B\boxplus B')
  \]
\end{Proposition}
\begin{proof}
  Let $\sX\to X$ be a $\m_p$-gerbe with cohomology class $\alpha$, and let $\sY\to Y$ be a $\m_p$-gerbe with cohomology class $\beta$. We regard $\sE$ as a $(-1,1)$-twisted sheaf on $\sX\times\sY$. As before, we consider the Cartesian diagram
\[
  \begin{tikzcd}
    \bP_{\sX\times\sY}(\sE)\arrow{r}{\pi'}\arrow{d}[swap]{g}&P\arrow{d}{f}\\
    \sX\times\sY\arrow{r}{\pi}&X\times Y
  \end{tikzcd}
\]
There is a universal quotient map $g^*\sE\twoheadrightarrow\cO(1)$, and $\cO(1)$ is $(-1,1)$-twisted. We set $\cO(p)=\pi'_*\cO(1)^{\otimes p}$. As in the proof of Proposition \ref{prop:integral1}, we may assume that $B$ and $B'$ have the property that $f^*(-B+B')=c_1^{\cry}(\cO(p))/p+h$ for some $h\in\H^2(P/W)$. By the corresponding result in the untwisted case (eg. \cite[Lemma 10.6]{Huy06}), it will suffice to show that $f^*(e^{-B\boxplus B'}\ch_{\sX\times\sY}(\sE))\in\H^*(P/W)$.
The same computation as in the proof of Proposition \ref{prop:integral1} shows that
\[
  f^*(e^{-B\boxplus B'}\ch_{\sX\times\sY}(\sE))=e^{-h}\ch(\pi'_*(g^*(\sE)\otimes\cO(-1)))
\]
The only primes dividing the denominators in right hand side of the above expression are $2$ and $3$, so this gives the result if $p\geq 5$.
\end{proof}
\section{Deformations of kernels of Fourier-Mukai transforms}\label{app:deformation}

In this appendix we record some results concerning deformations of kernels of Fourier-Mukai equivalences. The relevant deformation theory has been worked out by Reinecke in \cite{2017arXiv171100846R}, following \cite{LO15} in the non-twisted case. We apply this to complete the proof of Theorem \ref{thm:Mainprelim}, and to prove certain cases of Conjecture \ref{conj:mainPrelim2}.

Let $(X,\alpha)$, $(Y,\beta)$ be twisted supersingular K3 surfaces with B-field lifts $B,B'$, and suppose that $P^{\bullet}\in D(X\times Y,-\alpha\boxtimes\beta)$ is a perfect complex of twisted sheaves inducing a Fourier-Mukai equivalence
\[
  \Phi_{P^\bullet}\colon D(X,\alpha)\to D(Y,\beta)
\]
We will show that $P^{\bullet}$ deforms over a family of twisted surfaces connecting $(X\times Y,-\alpha\boxtimes\beta)$ to a product of twisted surfaces with trivial Brauer class. A special feature of the supersingular case is that there is a canonical choice of such a family. Specifically, we show in \cite[Lemma 2.2.5ff]{BL17} that if $X$ is a twisted supersingular K3 surface, then there exists a class 
\[
  \alpha_X\in\H^2(X\times\bA^1,\bG_m)
\]
with the property that for every $\alpha\in\Br(X)$, there exists a $k$-point $t\in\bA^1$ such that the restriction of $\alpha_X$ to $X\times t=X_{t}$ is equal to $\alpha$, and in particular the restriction of $\alpha_X$ to $X_0$ vanishes (there will be more than one such $t$ for a given $\alpha$).\footnote{Roughly speaking, the class $\alpha_X$ realizes $\bA^1$ as the moduli space of twisted supersingular K3 surfaces with underlying surface $X$ equipped with a certain level structure. This viewpoint is developed in detail in \cite{BL17}} In particular, all twisted supersingular K3 surfaces with underlying surface $X$ are deformation equivalent.

\begin{Proposition}\label{prop:deformations!!}
  Consider the universal class $\beta_Y\in\Br(Y\times\bA^1)$, and let $t\in\bA^1$ be a point such that the restriction of $\beta_Y$ to $Y_t$ is equal to $\beta$. There exists a connected \'{e}tale neighborhood $(U,u)\to (\bA^1,t)$ and
  \begin{enumerate}
      \item a family of K3 surfaces $X_U\to U$ such that $(X_U)_u\cong X$,
      \item a class $\alpha_U\in\Br(X_U)$ whose restriction to $(X_U)_u$ is equal to $\alpha$, and
      \item a kernel $\widetilde{P}^\bullet\in D(X_U\times Y_U,\alpha_U^{-1}\boxtimes\beta_U)$ inducing a Fourier-Mukai equivalence on each fiber, and such that $\widetilde{P}^\bullet_u\cong P^\bullet$.
  \end{enumerate}
\end{Proposition}
\begin{proof}
    Let $\widetilde{Y}\to Y\times\bA^1$ be a $\bG_m$-gerbe with cohomology class $\beta_Y$.
    Its underlying surface is the constant family $Y\times\bA^1\to\bA^1$, so determinant of any twisted sheaf on a fiber $\widetilde{Y}_x$ is unobstructed. Hence, the stack of simple perfect complexes $s\mathscr{T}w_{\widetilde{Y}/\bA^1}$ is smooth over $\bA^1$ (see for instance Definition 3.7 of \cite{2017arXiv171100846R}). The result then follows by the proof of Proposition 4.2 of \cite{2017arXiv171100846R}.
\end{proof}

To apply this, we record the following lemma. 
\begin{Lemma}\label{lem:itsalemma}
   Let $S$ be the spectrum of a DVR over $k$ with residue field $k$ and field of fractions $L$, $\sX\to S$ and $\sY\to S$ two families of supersingular K3 surfaces, $\alpha\in\Br(X)$ and $\beta\in\Br(Y)$ two Brauer classes, and $\sP^\bullet\in D(\sX\times \sY,\alpha^{-1}\boxtimes\beta)$ a perfect complex of twisted sheaves. The following are equivalent.
   \begin{enumerate}
       \item There exists a geometric point $t\in S$ and B-field lifts $B,B'$ of $\alpha_t$ and $\beta_t$ such that 
       \[
         v^{-B\boxplus B'}(\sP^{\bullet}_{t})\in\tN(\sX_{t}\times_{t}\sY_{t},-B\boxplus B')
       \] 
       \item For every geometric point $t\in S$ there exist B-field lifts $B,B'$ of $\alpha_t$ and $\beta_t$ such that
       \[
         v^{-B\boxplus B'}(\sP^{\bullet}_t)\in\tN(\sX_t\times \sY_t,-B\boxplus B')
       \]
   \end{enumerate}
\end{Lemma}
\begin{proof}
  The above notation was introduced in (\ref{eq:MoLatticesMoProblems}). Note that if either of the above is true for some choice of B-field lifts, it is in fact true for every choice of B-field lifts. 
  After possibly replacing $S$ with a connected \'{e}tale cover, we may find a complete lift $S'$ of $S$ over $W$ and elements $\sB\in\H^2(\sX/S'_K)$, $\sB'\in\H^2(\sY/S'_K)$ such that the restriction of $\sB$ to the special fiber is a B-field lift of $\alpha_0$, and the restriction of $\sB$ to any geometric generic fiber $\sX_{\overline{L}}$ is a B-field lift of $\alpha_{\overline{L}}$, and similarly, the restrictions of $\sB'$ give B-field lifts of the restrictions of $\beta$ in both special and geometric fibers. We consider the class
  \[
    \exp(\sB\boxplus-\sB') v^{-\sB\boxplus\sB'}(\sP^\bullet)\in\H^*(\sX\times_S\sY/S'_K)
  \]
  (the definition of the twisted Mukai vector makes sense for a twisted sheaf on a gerbe over any scheme). As $S$ is connected, this class is contained in the integral submodule $\H^*(\sX\times_S\sY/S')$ if and only if its restriction to any geometric point is. This gives the result.
  
\end{proof}

\begin{proof}[Proof of Theorem \ref{thm:Mainprelim}]
  It remains to show that the cohomological transform $\Phi^{\cry}_{v^{-B\boxplus B'}(P^{\bullet})}$ induces an isomorphism $\tH(X/W,B)\iso\tH(Y/W,B')$, and restricts to an isomorphism $\tN(X,B)\iso\tN(Y,B')$. To show this, it will suffice to prove that
  \[
    v^{-B\boxplus B'}(P^{\bullet})\in\tN(X\times Y,-B\boxplus B')
  \]
  As the square root of the Todd class of a K3 surface is integral, Proposition \ref{prop:integral2} already implies this when $p\geq 5$ (without using the assumptions that $X$ and $Y$ are supersingular, or that $\Phi_{P^{\bullet}}$ is an equivalence).
  
  We apply Proposition \ref{prop:deformations!!} to find a connected \'{e}tale neighborhood $(U,u)\to (\bA^1,t)$, a family of K3 surfaces $X_U\to U$, a class $\alpha_U\in\Br(X_U)$, and a kernel $\widetilde{P}^\bullet\in D(X_U\times Y_U,\alpha_U^{-1}\boxtimes\beta_U)$ extending our given data. Let $U'\to\bA^1$ be the normalization of $U$ in the function field of $\bA^1$. Let $R$ be the hensalization of $U'$ at a closed point mapping to $0\in\bA^1$. Let $L$ be the function field of $R$. Pulling back the data produced by Proposition \ref{prop:deformations!!}, we find a family of K3 surfaces $X_L$ over $L$, a class $\alpha_L\in\Br(X_L)$, and a kernel $P^\bullet_L\in D(X_L\times_L Y_L,\alpha_L^{-1}\boxtimes\beta_L)$, inducing a Fourier-Mukai equivalence on geometric fibers. In particular, this implies that $X_{\overline{L}}$ is supersingular. After replacing $R$ with a finite extension, we may therefore assume that $X_L$ is the generic fiber of a family $X_R$ of supersingular K3 surfaces over $R$ by \cite[Theorem 3]{RS83} (see also \cite[Theorem 5.2.1]{BL17} for $p=3$). Consider a flat extension $P^\bullet_R$ of $P^{\bullet}_L$ to $X_R\times_RY_R$. By Lemma \ref{lem:itsalemma}, we reduce to the case when $\beta=0$. Applying the same reasoning with the roles of $X$ and $Y$ swapped, we furthermore reduce to the case where $\alpha=0$. The result then follows from the well known fact that the Mukai vector of any complex on the product of two K3 surfaces is integral (see eg.  \cite[Lemma 10.6]{Huy06}). 
\end{proof}

We now discuss orientation. 
\begin{Lemma}\label{lem:this is a lemma}
   Let $S$ be the spectrum of a DVR with residue field $k$ and field of fractions $L$, $\sX\to S$ and $\sY\to S$ two families of K3 surfaces, and $P\in\tN(\sX\times \sY)_{\bQ}$ a class that induces an isomorphism $\Phi_{P_t}\colon\tN(\sX_t)_{\bQ}\iso\tN(\sY_t)_{\bQ}$ on each geometric fiber $r\in S$. The following are equivalent.
   \begin{enumerate}
       \item There exists a geometric point $t\in S$ such that $\Phi_{P_t}$ is orientation preserving.
       \item For every geometric point $t\in S$ the isometry $\Phi_{P_t}$ is orientation preserving.
   \end{enumerate}
\end{Lemma}
\begin{proof}
  As explained in \cite{2017arXiv171100846R}, an orientation on a lattice $N$ gives rise to a choice of connected component of the associated positive definite Grassmannian, and an isometry of lattices is orientation preserving if and only if it preserves the choice of connected components. The result follows.
\end{proof}

Using Lemma \ref{lem:this is a lemma} and arguments of \cite{MR2553878} and \cite{LO15}, we prove Conjecture \ref{conj:mainPrelim2} in certain special cases.
\begin{Theorem}\label{thm:nontwisted}
  Let $X$ and $Y$ be K3 surfaces over $k$, and suppose that $P^{\bullet}\in D(X\times Y)$ is a complex inducing a Fourier-Mukai equivalence
  \[
    \Phi_{P^{\bullet}}\colon D(X)\iso D(Y)
  \]
  The induced cohomological transform $\Phi^{\cry}_{v(P^{\bullet})}$ is orientation preserving.
\end{Theorem}
\begin{proof}
  By composing with standard derived equivalences (which are known to be orientation preserving), we reduce to the case where $\Phi_{P^\bullet}$ is filtered, in the sense of \cite{LO15}. The lifting results of \cite{LO15} combined with Lemma \ref{lem:this is a lemma} then reduce us to the case of K3 surfaces over the complex numbers, which follows from work of Huybrechts, Macr\'\i, and Stellari \cite{MR2553878}.
\end{proof}

\begin{Theorem}
  Suppose that $(X,\alpha)$ and $(Y,\beta)$ are twisted supersingular K3 surface such that $\sigma_0(X,\alpha),\sigma_0(Y,\beta)\leq 10$. If $P^{\bullet}\in D(X\times Y,\alpha^{-1}\boxtimes\beta)$ is a perfect complex of twisted sheaves that induces a Fourier-Mukai equivalence
  \[
    \Phi_{P^{\bullet}}\colon D(X,\alpha)\to D(Y,\beta)
  \]
  then for any B-field lifts $B,B'$ of $\alpha$ and $\beta$ the cohomological correspondence $\Phi^{\cry}_{v^{-B\boxplus B'}(P^{\bullet})}$ is orientation preserving.
\end{Theorem}
\begin{proof}
  Under our assumptions on the Artin invariant, we may find primitive isotropic vectors $v,w\in\tN(X,B)$ such that $v.w=-1$, and $v$ is effective. By results of \cite{BL17} (extending results of Mukai and Yoshioka over the complex numbers) the moduli space $X'=M_{(X,-\alpha)}(v)$ of $-\alpha$-twisted sheaves on $X$ with Mukai vector $v$ is a supersingular K3 surface, and the universal sheaf induces a derived equivalence $D(X,\alpha)\to D(X')$. As in \cite{HS04}, one shows that the induced cohomological correspondence is orientation preserving. Let $Y'$ be the same for $(Y,\beta)$. Now, given a Fourier-Mukai equivalence $D(X,\alpha)\to D(Y,\beta)$ that is not orientation preserving, we obtain by composition a Fourier-Mukai equivalence $D(X')\to D(Y')$ that is not orientation preserving, in contradiction to Theorem \ref{thm:nontwisted}.
\end{proof}

\begin{Remark}
  One can show that if $X$ is a K3 surface that is not supersingular and $\alpha\in\Br(X)$ is any Brauer class, then $(X,\alpha)$ lifts to characteristic 0. The argument of Theorem \ref{thm:nontwisted} then extends to prove Conjecture \ref{conj:mainPrelim2} for all twisted K3 surfaces in positive characteristic that are not supersingular. Thus, the only case of Conjecture \ref{conj:mainPrelim2} left open is that of non-trivial Brauer classes on supersingular K3 surfaces of Artin invariant $10$.
\end{Remark}


\end{document}